\documentclass[11pt,leqno,twoside]{article}

\usepackage{amsfonts,amsmath,amsthm,amssymb}
\usepackage{enumerate}
\usepackage{graphics,graphicx,subfigure}

\usepackage{color}
\usepackage{todonotes}
\usepackage{cancel}
\usepackage{url}
\usepackage{hyperref}
\usepackage{makeidx}
\usepackage{showidx}
\usepackage{multicol}        
\usepackage{xspace}
\usepackage{stmaryrd}        
\usepackage{pifont}          
\usepackage{fancybox}        
\usepackage{bm}

 \usepackage{fancyhdr}
\theoremstyle{plain}
 \theoremstyle{definition}
 \newtheorem{lem}{Lemma}
 \newtheorem{defn}[lem]{Definition}
 \newtheorem{thm}[lem]{Theorem}
 \newtheorem{prop}[lem]{Proposition}
 \newtheorem{cor}[lem]{Corollary}
 \newtheorem{notn}[lem]{Notations}
 \newtheorem{pb}[lem]{Problem}
 \newtheorem{form}[lem]{Formulae}
 
 \newtheorem*{rk}{Remark}
 \newtheorem*{com}{Comment}
 \newtheorem*{ex}{Example}
 \theoremstyle{remark}

 \newcommand{\blem}{\begin{lem}}
 \newcommand{\elem}{\end{lem}}
 \newcommand{\bdefn}{\begin{defn}}
 \newcommand{\edefn}{\end{defn}}
 \newcommand{\bthm}{\begin{thm} }
 \newcommand{\ethm}{\end{thm}}
 \newcommand{\bprop}{\begin{prop}}
 \newcommand{\eprop}{\end{prop}}
 \newcommand{\bcor}{\begin{cor}}
 \newcommand{\ecor}{\end{cor}}
 \newcommand{\bnotn}{\begin{notn}}
 \newcommand{\enotn}{\end{notn}}
 \newcommand{\bpb}{\begin{pb}}
 \newcommand{\epb}{\end{pb}}
 \newcommand{\bform}{\begin{form}}
 \newcommand{\eform}{\end{form}}
 \newcommand{\brk}{\begin{rk}}
 \newcommand{\erk}{\end{rk}}
 \newcommand{\bcom}{\begin{com}}
 \newcommand{\ecom}{\end{com}}
 \newcommand{\bex}{\begin{ex}}
 \newcommand{\eex}{\end{ex}}
 \newcommand{\bpf}{\begin{proof}}
 \newcommand{\epf}{\end{proof}}

\newcommand{\cD}{\mathcal{D}}

\newcommand{\cK}{\mathcal{K}}

\newcommand{\cN}{\mathcal{N}}

\newcommand{\cV}{\mathcal{V}}

\newcommand{\cX}{\mathcal{X}}

\newcommand{\bC}{\mathbb{C}}

\newcommand{\bK}{\mathbb{K}}

\newcommand{\bR}{\mathbb{R}}

\newcommand{\be}{\begin{equation}}
\newcommand{\ee}{\end{equation}}
\newcommand{\bal}{\begin{align}}
\newcommand{\eal}{\end{align}}
\newcommand{\ba}{\begin{align*}}
\newcommand{\ea}{\end{align*}}
\newcommand{\bmx}{\begin{matrix}}
\newcommand{\emx}{\end{matrix}}
\newcommand{\bbmx}{\begin{bmatrix}}
\newcommand{\ebmx}{\end{bmatrix}}
\newcommand{\bpmx}{\begin{pmatrix}}
\newcommand{\epmx}{\end{pmatrix}}
\newcommand{\bvmx}{\begin{vmatrix}}
\newcommand{\evmx}{\end{vmatrix}}

\newcommand{\ol}{\overline}
\newcommand{\wh}{\widehat}
\newcommand{\wt}{\widetilde}
\newcommand{\f}{\frac}
\newcommand{\df}{\dfrac}

\newcommand{\imp}{\Longrightarrow}

\newcommand{\inc}{\subseteq}

\newcommand{\Id}{\mathrm{Id}}

\newcommand{\argmin}{{\rm argmin}\,}

\newcommand{\minimize}[1]{\underset{#1}{\rm minimize}\,}

\newcommand{\la}{\lambda}
\newcommand{\La}{\Lambda}
\newcommand{\eps}{\varepsilon}
  
\title{\vspace{-11mm}Worst-Case Learning under a  Multi-fidelity Model
\medskip\hrule height 1.2pt \vspace{-6mm}}
\author{Simon Foucart\footnote{Texas A\&M University} \, and Nicolas Hengartner\footnote{Los Alamos National Laboratory\\S. F. is partially supported by grants from the NSF (DMS-2053172) and from the ONR (N00014-20-1-2787).
S. F. is also grateful for a Development Fellowship from the Texas A\&M University System National Laboratories Office, which made visits to Los Alamos National Laboratory possible during 2023--24.}}
\date{\vspace{-6mm}\rule{100mm}{0.8pt}}

\newcommand\shorttitle{Worst-Case Learning under a  multi-fidelity Model}
\newcommand\authors{S. Foucart, N. Hengartner}

\begin{document}
\maketitle

\vspace{-15mm}
\begin{abstract}
Inspired by multi-fidelity methods in computer simulations,
this article introduces procedures to design surrogates for the input/output relationship of a high-fidelity code.
These surrogates should be learned from runs of both the high-fidelity and low-fidelity codes
and be accompanied by error guarantees that are deterministic rather than stochastic.
For this purpose, the article advocates 
a framework tied to a theory focusing on worst-case guarantees, namely Optimal Recovery.
The multi-fidelity considerations triggered 
new theoretical results in three scenarios:
the globally optimal estimation of linear functionals,
the globally optimal approximation of arbitrary quantities of interest in Hilbert spaces,
and their locally optimal approximation, still within Hilbert spaces.
The latter scenario boils down to the determination of the Chebyshev center for the intersection of two hyperellipsoids.
It is worth noting that the mathematical framework presented here, together with its possible extension,
seems to be relevant in several other contexts briefly discussed.
\end{abstract}



\vspace{-5mm}
\begin{center}
\rule{100mm}{0.8pt}
\end{center}

\section{Introduction}
High-fidelity computer models of complex  
physical systems provide some unprecedentedly detailed insights into processes ranging 
from climate response to anthropogenic greenhouse gas emissions~\cite{edwards2011history},  
to the dynamics of turbulence and combustion of scramjets~\cite{hoste2017numerical} and rocket motors~\cite{chen2011multiphysics}, 
to modeling the dynamic responses of materials to shocks~\cite{steinhauser2009review}.  
These models couple elementary
physical, chemical, and biological processes at scale, leading to  nonlinear dynamics that often reveal emergent phenomena~\cite{crutchfield1994calculi}.    

The computational cost of running high-fidelity models can be prohibitively high, however,
which may limit their usefulness for some applications. 
Low-fidelity models of the high-fidelity codes rely on streamlined assumptions such as linearization, lower-dimensional approximations, simplified physics, complexity management limiting the number and order of interactions, and domain coarsening,  
to significantly reduce the computational burden.  
{\em Multi-fidelity} models \cite{fernandez2016review, peherstorfer.2018} combine the results from both low- and high-fidelity models,  
using the low-fidelity ones to extensively explore an approximation of the relationship between inputs and outputs
and capitalizing on the higher fidelity ones to finetune this relationship.
We note that the notions of high-fidelity and low-fidelity are only defined relative to one another.

{\em Surrogate models} are data-driven approximations of the relationship between inputs and outputs.  
They are essential ingredients in optimization \cite{koziel.2011} and in estimation and uncertainty quantification of computer models \cite{queipo.2005,santnerBook}.
Developing surrogate models for multi-fidelity computer models
is challenging and remains an active field of research
\cite{Forrester.2007,Park.2017}.
While modern machine learning methods \cite{Zhang.2021,meng.2020,Song.2019} have successfully combined low- and high-fidelity data, 
they are not well equipped to bound the pointwise prediction error.  
Multi-fidelity models based on Gaussian processes \cite{SongJ.2019, haaland2018} can provide stochastic bounds on the error.  
Yet, for some applications, deterministic bounds are desirable.
The theory of Optimal Recovery---the centerpiece of this article, going back to the seventies \cite{MicRiv}---is tailored to produce such deterministic bounds,
e.g. for surrogates of computer simulations.  
These bounds offer a complementary view of the uncertainty associated 
with the surrogate, with a different interpretation but often similar ranges.  
See \cite{frenklach.2016} for an example of comparison between stochastic and deterministic uncertainty quantification.  

This article contributes to the existing body of research on surrogate modeling by extending the theory of Optimal Recovery to the design of surrogates
for multi-fidelity computer simulations.  
The benefit of this approach is the explicit determination of upper bounds for the error of prediction,
or even of full approximation, of the true input/output relationship by a surrogate constructed precisely to make this error small under some prior assumption about the relationship.
Thus, it provides novel deterministic uncertainty quantification tools for multi-fidelity computer models.

In mathematical terms,
let the input/output relationship for the high-fidelity code be represented by a function $f_0$
and let a function $f_1$ represent this relationship for the low-fidelity code.
The prior assumption takes the form of a condition $f_0 \in \cK_0$,
together with a condition $f_0 - f_1 \in \cK_1$
expressing the belief that, in a multi-fidelity setting, the bias $f_0 - f_1$ should be well-behaved, hopefully of lower complexity than $f_0$.
Executing the high- and low-fidelity codes provides the data  
$$
y_{0,i} = \lambda_{0,i}(f_0), 
\quad i=1,\ldots,m_0,
\qquad \mbox{and} \qquad
y_{1,i} = \lambda_{1,i}(f_1), 
\quad i=1,\ldots,m_1,
$$
where the linear functionals $\la_{0,i}, \la_{1,i}$ are often (but not always) point evaluations.  
These~data serve as constraints on the unknown functions, expressed via
the sets $\cD_0= \{ f: \lambda_{0,i}(f)=y_{0,i} \mbox{ for all }i\}$
and $\cD_1 = \{ f : \lambda_{1,i}(f)=y_{1,i} \mbox{ for all }i\}$.
The error of full approximation, say,
is obtained as the optimal value of the 
min-max program
\begin{equation}\label{eq:min.max.A}
\inf_h \; \sup \{ \|f_0-h\| : f_0 \in \cK_0 \cap \cD_0, f_1-f_0 \in \cK_1, f_1 \in \cD_1 \},
\end{equation}
while the minimizer $h$ amounts to the optimal surrogate.

This formalism extends to more than two fidelity levels.  
Although conceptually straightforward, 
such an extension introduces theoretical challenges 
for solving suitable generalizations of~(\ref{eq:min.max.A}),
so our results primarily focus on two levels.
Nonetheless, we want to make the point that the core mathematical problem arises in 
contexts beyond the multi-fidelity setting considered 
in this paper, 
and as such is of interest to a broad audience.  
Below is a non-exhaustive list of four examples that can benefit from our Optimal Recovery perspective.

\vspace{-2mm}
\textit{Accelerated life testing.}  
Accelerated life testing subjects items to conditions that enhances the physical and chemical processes associated with aging.  
For example, increasing the temperature of explosives promotes chemical degradation. 
We can quantify this effect in many ways, 
e.g. gas chromatography to characterize volatile and semi-volatile mixtures \cite{freye2023a}, possibly complemented with time-of-flight mass spectroscopy \cite{freye2023b}, and `dynamic testing' to measure the energetic release.  
These destructive experiments capture physical and chemical characteristics that are differentially impacted by aging.  
For each experiment, we are interested in learning the expected response $f_t(x)$ to the $t$-th experiment, 
where $x$ represents aging.
The framework for Optimal Recovery in multi-fidelity experiments applies here, 
allowing us to combine and learn from multiple related tests.

\vspace{-2mm}
{\em Multi-physics sensing}.  
Nuclear particles, such as electrons, muons, neutrons and protons  penetrate and interact with matter in different but related ways. 
For example, muons  and electrons interact with matter through electro-magnetic forces (Coulomb scattering) \cite{schultz2007, cuellar2009},
whereas neutrons scatter when colliding with the nucleus of atoms.  
Both of these processes are enhanced in high-Z materials and there is interest in combining the probing ability of multiple types of particles.  
The presented Optimal Recovery framework applies here, too, as we seek to learn the function $f_t(x,z)$, where $x$ is the spatial position in a probed object
and $z$ are characteristics of the individual particles of type~$t$.

\vspace{-2mm}
\textit{Panel data.}  
Panel data consists of short time series observed on multiple individuals, with common dynamic for its temporal evolution
\cite{hsiao.2022}, possibly depending on covariates $x$.  For times $t =0, 1,\ldots, T$,  we denote by $f_t(x)$ the response at time $t$ as a function of covariate $x$.  
The temporal evolution of the process makes it plausible that measurements at time $t$ may inform, or at least bound, the measurements at time $t'$.  
That is, we can again leverage the theory in this paper to combine data from multiple time points to estimate the dependence of the response at time $t$ on the covariate $x$, even if the measurements are from different items. 

\vspace{-2mm}
\textit{Theoretical tool for transfer learning.}  Transfer learning has emerged as a powerful machine-learning paradigm to combine related yet distinct datasets,  
see the recent review articles \cite{weiss.2016, zhuang.2020}.
The formalism in this paper provides a way to quantify the benefit of combining datasets by computing uncertainty bound of surrogates with and without the addition of `related datasets'.

\section{Preliminaries, Problem Formulation, and Summary of Results} 
\label{SecForma}

In this section, we fix  the notation, introduce the concepts, and provide a self-contained overview of {\em Optimal Recovery},
specialized to the problem of building surrogates in our
multi-fidelity context.  
In addition to \cite{MicRiv},
the books \cite{NovWoz,BookDS}  can provide the  interested reader with a more modern discussion about the theory and applications of Optimal Recovery.

In this field, the goal is to exploit the assumption that an unknown  target function 
$f_0$ belongs to a given subset $\cK \inc F$ of a Banach space
$F$, together with observed constraints $y_{0,1}, \ldots, y_{0,m_0}$, which we call data, to bound a quantify of interest $Q(f_0)$.  
For simplicity, the latter is taken to be a linear map $Q:F \to Z$ into another Banach space $Z$.  
This framework includes the problem of full approximation of the function $f_0$ by setting $Q=\Id$, the identify map.   
The Banach space $F$ already encapsulates some {\em a priori} assumptions  
about $f_0$ and the so-called model set $\cK_0$ further reflects our educated knowledge about realistic objects.  

Throughout this paper, we will suppose that the data 
constraints 
$$
y_{0,i} = \la_{0,i}(f_0),
\qquad i=1,\ldots,m_0,
$$
are obtained by applying linear functionals $\la_{0,1},\ldots, \la_{0,m_0}$ to $f_0$.  
It is convenient to vectorize these data and use the notation  $y_0 = \La_0 f_0 = [\la_{0,1}(f_0);\ldots;\la_{0,m_0}(f_0)] \in \bR^{m_0}$ and refer to $\La_0$ as the observation map from
$F$ to $\bR^{m_0}$.  
A special case of particular interest is when
$\lambda_{0,i}(f_0)=f_0(x^{(0,i)})$ is a point evaluation functional.  
This evokes similarities with Statistical Learning Theory.
However, the main difference is that the $x^{(i)}$'s are regarded as fixed here and not as independent realizations of a random variable.
As a result, the performance of a learning/recovery procedure cannot be assessed in an average case.
Thus, one opts for an assessment focusing on the worst case, see \eqref{LWCEOri}~and~\eqref{GWCEOri}
below, which is a distinctive feature of the Optimal Recovery framework. 

In our specific setting,
let us recall that multi-fidelity supplies an extra dimension via some side information on an object $f_1$
related to $f_0$.
Indeed, the latter is available through its own {\em a priori} assumption on the bias $f_0 - f_1$ expressed via
$$
f_0 - f_1 \in \cK_1.
$$
For instance, if we have knowledge that $f_1$ is $\eps$-close to $f_0$,
then $\cK_1$ would simply be a ball of radius~$\eps$.
As for the data constraints,
they are observations made directly on $f_1$,
i.e.,
$$
y_{1,i} = \la_{1,i}(f_1),
\qquad
i=1,\ldots,m_1.
$$
These data are again vectorized as $y_1 = \La_1 f_1$,
where $\La_1: F \to \bR^{m_1}$ is a linear map.

To estimate $Q(f_0) \in Z$,
we want to produce an approximant $\wh{z}$ built from the combined data $y = [y_0; y_1] \in \bR^{m_0+m_1}$.
In other words, we want to construct a mapping $\Delta$ from $\bR^{m_0+m_1}$ to $Z$,
which we call a recovery map.
Note that it may be cognizant of the model sets $\cK_0$ and $\cK_1$.
We are interested in {\em optimal} recovery maps,
the meaning of which being tied to how the performance of a generic $\Delta$ is assessed.
The value of $\|Q(f_0) - \Delta(y) \|_Z$ evidently quantifies the recovery error for fixed $f_0$ and $f_1$,
but since $f_0$ and $f_1$ are unknown,
we settle on a worst-case viewpoint to quantify performance by taking the supremum over all $f_0$ and $f_1$ that are consistent with model and data.
This leads to the following two notions of worst-case recovery error: \vspace{-5mm}
\begin{itemize}
\item the local worst-case error, at a fixed $y=[y_0;y_1] \in \bR^{m_0+m_1}$, is defined as
\be
\label{LWCEOri}
{\rm lwce}_y(z) = \sup_{\substack{f_0 \in \cK_0, f_0 - f_1 \in \cK_1 \\ \La_0 f_0 = y_0, \, \La_1 f_1 = y_1}} \|Q(f_0) - z \|_Z;
\ee
\item the global worst-case error is defined as
\be
\label{GWCEOri}
{\rm gwce}(\Delta) = \sup_{f_0 \in \cK_0, f_0 - f_1 \in \cK_1} \|Q(f_0) - \Delta([\La_0 f_0;\La_1 f_1]) \|_Z.
\ee
\end{itemize}
A globally, resp. locally, optimal recovery map $\Delta^{\rm opt}: \bR^{m_0+m_1} \to Z$ is a map that minimizes ${\rm gwce}(\Delta)$ over all $\Delta: \bR^{m_0+m_1} \to Z$,
resp. such that $\Delta^{\rm opt}(y)$ minimizes 
${\rm lwce}_y(z)$ over all $z \in Z$ at every $y \in \bR^{m_0+m_1}$.
It is straightforward to see that a locally optimal recovery map is automatically globally optimal,
implying that globally optimal recovery maps are  somewhat easier to come by than locally optimal recovery maps.
This will also be exemplified by the representative cases investigated in the remainder of this article.
Treated by increasing order of difficulty,
they leverage some techniques that have recently proved useful in Optimal Recovery.
Below is an summary of some results for the three scenarios we considered.

In the first scenario (Section \ref{SecLF}),
the quantity of interest is required to be a linear functional.
For instance, in the space $F=C(\cX)$ of continuous functions on a compact set $\cX$,
suppose that the observations $\la_{0,i}(f_0) = f_0(x^{(0,i)})$ and $\la_{1,i}(f_1) = f_1(x^{(1,i)})$
are point evaluations
and that the quantity of interest $Q(f_0) = f_0(x)$ is the evaluation at another point $x \in \cX$,
i.e., we are trying to predict $f_0$ at new point.
Suppose also that the model sets are based on approximation capabilities by linear subspaces $\cV_0,\cV_1$ of $C(\cX)$ with paramaters $\eps_0,\eps_1 > 0$,
i.e.,
$$ 
\cK_0 = \{ g \in C(\cX): {\rm dist}(g,\cV_0) \le \eps_0 \}
\qquad \mbox{and} \qquad
\cK_1 = \{ g \in C(\cX): {\rm dist}(g,\cV_1) \le \eps_1 \}.
$$
Then,
a {\em globally} optimal estimation map takes the form
$$
y = [y_0; y_1] \in \bR^{m_0 + m_1} \mapsto 
\sum_{i=1}^{m_0} a^{\rm opt}_{0,i} y_{0,i}
+ \sum_{i=1}^{m_1} a^{\rm opt}_{1,i} y_{1,i} \in \bR,
$$
where the vector $a^{\rm opt} = [a^{\rm opt}_0;a^{\rm opt}_1] \in \bR^{m_0} \times \bR^{m_1}$ is solution to the $\ell_1$-optimization program
$$
\minimize{[a_0; a_1] \in \bR^{m_0+m_1}} 
\big[
\eps_0  \, \|a_0\|_1  + (\eps_0 + \eps_1) \, \|a_1\|_1 \big]
\qquad \mbox{s.to }
M^{(0)} [a_0; a_1] = b
\; \mbox{ and } \; M^{(1)} a_1 = 0
$$
for some explicit matrices $M^{(0)}, M^{(1)}$, and vector $b$.
This result is in fact an instantiation of a more general result (Corollary \ref{CorEstLFinCX}).
Indeed, in the case of estimation of a linear functional $Q$,
the result can be extended to any number of model sets $\cK_0,\cK_1,\ldots,\cK_T$---which do not need to be approximability sets either---and therefore apply to the four examples mentioned in the introduction.
This extension can also be viewed as a first step towards a `dynamical Optimal Recovery',
where the object of interest is some time-varying $f_t$ initially modeled by $f_0 \in \cK_0$
and whose evolution is modeled by $f_{t-1}-f_t \in \cK_t$, $t \ge 1$.
As a stylized example, we are thinking of the estimation of average temperatures at a past time $t=0$ given observations at this time $t=0$ and at subsequent times $t \ge 1$.
For this climate application,
the use of Optimal Recovery in the `static' framework was proposed in \cite{FHMPW}.

In the second scenario (Section \ref{SecGlo}), 
there is no restriction on the quantity of interest $Q: F \to~Z$
(except that it is a linear map),
but there are restrictions on the spaces $F$ and $Z$,
namely, they should be Hilbert spaces.
Thus, when the objects $f_0$ and $f_1$ are functions observed via point evaluations,  reproducing kernel Hilbert spaces provide the right framework.
For arbitrary linear observation maps $\La_0: F \to \bR^{m_0}$ and $\La_1: F \to \bR^{m_1}$,
suppose that the model sets are given by 
$$
\cK_0  = \{ g \in F: \|P_0 g \| \le \eps_0 \}
\qquad  \mbox{and} \qquad
\cK_1  = \{ g \in F: \|P_1 g \| \le \eps_1 \},
$$
where the linear operators $P_0,P_1$ map into Hilbert spaces, too.
Then,
a {\em globally} optimal recovery map is produced 
by constrained regularization with an explicitly determined parameter.
Precisely, it is obtained as $Q \circ \Delta_{\tau^\sharp}: \bR^{m_0 + m_1} \to Z$,
where $\Delta_\tau$ is defined for $\tau \in [0,1]$ by $\Delta_\tau([y_0;y_1]) = f^\tau_0 \in F$ and 
$$
[f^\tau_0;f^\tau_1]
= \underset{f = [f_0; f_1] \in F \times F}{\argmin \;}
\left[ (1-\tau) \|P_0 f_0\|^2 + \tau \|P_1(f_0 - f_1) \|^2 \right]
\qquad \mbox{s.to } \La_0 f_0 = y_0, \; \La_1 f_1 = y_1.
$$
As for the parameter $\tau^\sharp \in [0,1]$,
it is selected as $\tau^\sharp = c_1^\sharp/ (c_0^\sharp + c_1^\sharp)$,
where $c_0^\sharp , c_1^\sharp \ge 0$ are solutions to the semidefinite program
\begin{align*}
\minimize{c_0,c_1 \ge 0} &  \qquad c_0 \eps_0^2 + c_1 \eps_1^2\\
 \mbox{s.to } & \qquad  c_0 \|P_0 f_0\|^2 + c_1 \|P_1(f_0 - f_1) \|^2 \ge \|Qf_0\|^2
\quad \mbox{for all } f_0 \in \ker(\La_0),
f_1 \in \ker(\La_1).
\end{align*}
This statement, exploiting a recent result from \cite{L1Err},
does not extend to any number of model sets.
This is because the underlying tool---Polyak's S-procedure, see \cite{Pol}---is similarly limited in the number of quadratic constraints it can handle.
Beyond the multi-fidelity setting,
situations that seem pertinent to this scenario
involve destructive observations.
As a stylized example, we can consider $f_0$ as a graph signal,
i.e., a function defined on the finitely many vertices of a graph,
whence it is relevant to take the operator $P_0$ describing $\cK_0$ as the square-root of the graph Laplacian.
We can also imagine that the very fact of observing $f_0$ is destructive, in the sense that it alters $f_0$ into a closeby $f_1$,
whence it is relevant to take the operator $P_1$ as the identity.
In a non-destructive framework,
the graph-signal application of Optimal Recovery was proposed in \cite{FouLiaVel}.

In the third scenario (Section \ref{SecLoc}), the setting is similar to the one from the second scenario,
except that a {\em locally} optimal recovery map is targeted instead of a globally optimal one.
In the spirit of~\cite{FouLia-CAMDA},
which extended the result from \cite{BecEld} to an arbitrary linear quantity of interest $Q \not= \Id$,
we first show that the minimal local worst-case error 
(aka local radius of information or Chebyshev radius)
is upper bounded by the square root of the optimal value of the semidefinite program
\begin{align}
\label{SDPIntro}
 \minimize{b,c_0,c_1 \ge 0} & \; 
c_0 \big( \eps_0^2 - \|P_0 \La_0^\dagger y_0\|^2 \big) + c_1 \big( \eps_1^2 - \|P_1(\La_0^\dagger y_0 - \La_1^\dagger y_1)\|^2 \big) + b\\
\nonumber
\mbox{s.to } &
\bbmx
c_0 M_0 + c_1 M_1 & \vline & - c_1 M_1 \\
\hline
- c_1 M_1 & \vline & c_1 M_1
\ebmx \succeq 
\bbmx
N_0 & \vline & 0\\
\hline 0 & \vline & 0 
\ebmx
\\
\nonumber
\mbox{and } &
\bbmx
c_0 M_0 + c_1 M_1 & \vline & - c_1 M_1
& \vline & 
c_0 L_0 y_0 + c_1 L_1 y
\\
\hline
- c_1 M_1 & \vline & c_1 M_1 & \vline & - c_1 L_1 y\\
\hline 
(c_0 L_0 y_0 + c_1 L_1 y)^* & \vline & - (c_1 L_1 y)^* & \vline & b
\ebmx \succeq 
0,
\end{align}
where $L_0, L_1,M_0, M_1$ are expressed in terms of $P_0, P_1,\La_0, \La_1$, 
while $N_0$ is expressed in terms of $Q,\La_0$.
Next, we show that this upper bound can actually be equal to the genuine minimal local worst-case error under some additional assumptions, e.g. $\ker(\La_0) \inc \ker(P_1)$ or $\ker(\La_0) \inc \ker(\La_1)$.
The locally optimal recovery map is still produced by constrained regularization:
it takes the form $Q \circ \Delta_{\tau^y}$,
where the parameter $\tau^y \in [0,1]$ now depends on $y$ via $\tau^y = c_1^y / (c_0^y + c_1^y)$,
with $b^y,c_0^y,c_1^y \ge 0$ being solutions to \eqref{SDPIntro}.
Beyond the multi-fidelity setting,
we can envisage a stylized example in which $f_0$ represents a physical entity and $f_1$ represents its digital twin.
The observation functionals $\la_{1,i}$'s are chosen among the $\la_{0,i}$'s---hence ensuring that $\ker(\La_0) \inc \ker (\La_1)$---for the purpose of confirming that $f_1$ is a convincing twin of $f_0$. 

The third scenario directly relies on more general results about Chebyshev radii and centers for a model set given by the intersection of two hyperellipsoids.
We isolate these seemingly novel results in Section \ref{SecCCtwoHyp}.
Not only do they supply alternative---and arguably simpler (compare to~\cite{FouLia})---arguments for the exact determination of other Chebyshev centers, 
but they also provide a {\em necessary and sufficient} orthogonality condition for a semidefinite-relaxation upper bound to agree with the true Chebyshev radius.
We refer to Theorem \ref{ThmCCtwoHyp} for the precise statement.

\section{Globally Optimal Estimation of Linear Functionals}
\label{SecLF}

This section deals with the first scenario,
where the quantity of interest $Q$ is a linear functional defined on an arbitrary Banach space $F$.
Exploiting the classical theory of Optimal Recovery,
we can swiftly guarantee that, under reasonable assumptions,
there is a globally optimal recovery map which is linear.
Our main focus will be on the efficient construction of such a map.
We shall do so in an extended framework
where more than two model sets can be involved.
Namely,
the object $f_0$ comes with related objects $f_1,\ldots,f_T$,
which altogether are available through the model assumptions
$$
f_0 \in \cK_0,
\qquad 
f_{t-1} - f_{t} \in \cK_t, 
\quad t=1,\ldots,T, 
$$
as well as the observed data
$$
y_0 = \La_0 f_0 \in \bR^{m_0},
\qquad 
y_t = \La_t f_t \in \bR^{m_t},
\quad t=1,\ldots,T.
$$
To simplify our analysis,
it is useful to reduce this extended framework to a more familiar one by considering
the compound vector 
$$
f = \bbmx f_0 \\ \hline f_1 \\ \hline \vdots \\ \hline f_T \ebmx \in F^{T+1}.
$$
The model assumptions on this vector simply become $f \in \cK$, where
$$
\cK = \left\{
\bbmx f_0 \\ \hline f_1 \\ \hline \vdots \\ \hline f_T \ebmx \in F^{T+1} : f_0 \in \cK_0, f_{t-1} - f_t \in \cK_t, \; t=1,\ldots,T
\right\},
$$
while the observed data takes the condensed form,
with $m:=m_0 + m_1 + \cdots + m_T$,
$$
y = \La f \in \bR^{m},
\qquad \mbox{where} \quad 
y = \bbmx y_0 \\ \hline y_1 \\ \hline \vdots \\ \hline y_T \ebmx 
\quad \mbox{and} \quad
\La  = \bbmx
\La_0 & \vline & 0 & \vline & \cdots & \vline & 0\\
\hline
0 & \vline & \La_1 & \vline &  & \vline & \vdots\\
\hline
\vdots & \vline & & \vline & \ddots & \vline & 0 \\
\hline
0 & \vline & \cdots & \vline & 0 & \vline & \La_T 
\ebmx.
$$
This notation allows us to write the worst-case errors of a recovery map $\Delta: \bR^{m} \to \bR$ in the simpler form
\be
\label{ReformGWCE_Q}
{\rm gwce}(\Delta)  = \sup_{f \in \cK} \big| \wt{Q}(f) - \Delta(\La f) \big|,
\ee
where $\wt{Q}: F^{T+1} \to \bR$ is the linear functional defined by $\wt{Q}([f_0;f_1;\ldots;f_T]) = Q(f_0)$.
This reduction facilitates the derivation of the following abstract statement,
soon to be specialized in more tangible situations.

\bthm
\label{ThmLFAbs}
Suppose that the model sets $\cK_0,\cK_1,\ldots,\cK_T$ are symmetric and convex.
If the quantity of interest $Q$ is a linear functional,
then a globally optimal recovery map is given by
$$
\Delta^{\rm opt}: y = [y_0; y_1; \ldots; y_T] \in \bR^{m} \mapsto 
\sum_{i=1}^{m_0} a^{\rm opt}_{0,i} y_{0,i}
+ \sum_{t=1}^T \sum_{i=1}^{m_t} a^{\rm opt}_{t,i} y_{t,i} \in \bR,
$$
where the vector $a^{\rm opt} = [a^{\rm opt}_0; a^{\rm opt}_1; \ldots; a^{\rm opt}_T] \in \bR^{m}$ is a solution to
\begin{align}
\label{MiniAbs}
\minimize{a = [a_0; a_1; \ldots; a_T] \in \bR^{m}} 
\sup_{g_0 \in \cK_0}
\left| 
\left( Q - \sum_{t=0}^T \sum_{i=1}^{m_t} a_{t,i} \la_{t,i}\right)(g_0) \right|  + \sum_{s=1}^T \sup_{g_s \in \cK_s}
\left| \left( \sum_{t=s}^T \sum_{i=1}^{m_t} a_{t,i} \la_{t,i} \right)(g_s)
\right|.
\end{align}
\ethm

\bpf
Given that $\cK_0,\cK_1,\cdots,\cK_T \inc F$ are all symmetric and convex, so is $\cK \inc F^{T+1}$.
Moreover, $Q$ being a linear functional, 
$\wt{Q}$ is also a linear functional.
Thus, according to the reformulation~\eqref{ReformGWCE_Q} of the global worst-case error, 
the foundational result of Smolyak
(see e.g. \cite[Theorem~4.7]{NovWoz} or \cite[Theorem~9.3]{BookDS})
ensures that,
among the recovery maps miminizing ${\rm gwce}(\Delta)$,
there exists one which is linear,
say of the form $\Delta_a = \langle a, \cdot \rangle$,
i.e.,
$\Delta_a([y_0;y_1;\ldots;y_T]) = \sum_{i=1}^{m_0} a_{0,i} y_{0,i} + \sum_{t=1}^T \sum_{i=1}^{m_t} a_{t,i} y_{t,i}$.
To find a globally optimal recovery map,
it is therefore enough to minimize ${\rm gwce}(\Delta_a)$ over all $a = [a_0; a_1; \ldots; a_T] \in \bR^{m}$.
Coming back to the original formulation of the global worst-case error and substituting the specific form of $\Delta_a$, we obtain
\begin{align*}
{\rm gwce}(\Delta_a) & = \sup_{\substack{f_0 \in \cK_0\\ f_{t-1} - f_t \in \cK_t, t = 1,\ldots,T}}
\left| 
Q(f_0) - \left(
 \sum_{i=1}^{m_0} a_{0,i} \la_{0,i}(f_0) + \sum_{t=1}^T \sum_{i=1}^{m_t} a_{t,i} \la_{t,i}(f_t) 
\right)
\right|\\
& = 
\sup_{\substack{f_0 \in \cK_0\\ f_{t-1} - f_t \in \cK_t, t=1,\ldots,T}}
\left| 
\left( Q -  \sum_{t=0}^T\sum_{i=1}^{m_t} a_{t,i} \la_{t,i}\right)(f_0) 
+ \left( \sum_{t=1}^T\sum_{i=1}^{m_t} a_{t,i} \la_{t,i} \right)(f_0-f_t)
\right| .
\end{align*}
Introducing temporarily $\mu := Q -  \sum_{t=0}^T\sum_{i=1}^{m_t} a_{t,i} \la_{t,i}$
and $\nu_t := \sum_{i=1}^{m_t} a_{t,i} \la_{t,i}$ for $t=1,\ldots,T$,
while writing $g_0 = f_0$ and $g_s = f_{s-1}-f_s$ for $s=1,\ldots,T$,
the quantity inside the absolute values is
$$
\mu(g_0) + \sum_{t=1}^T \nu_t \left( \sum_{s=1}^t g_s \right)
 = \mu(g_0) + \sum_{s=1}^T \left( \sum_{t=s}^T \nu_t \right)(g_s).
$$
Therefore, we arrive at 
$$
{\rm gwce}(\Delta_a)  =
\sup_{\substack{g_0 \in \cK_0\\ g_s \in \cK_s, s=1,\ldots,T}}
\left| \mu(g_0) + \sum_{s=1}^T \left( \sum_{t=s}^T \nu_t \right)(g_s) \right|
= \sup_{g_0 \in \cK_0}
\left| \mu(g_0)  \right|
+ \sum_{s=1}^T \sup_{ g_s \in \cK_s}
\left| \left( \sum_{t=s}^T \nu_t \right)(g_s) \right|,
$$
where the last step relied on the symmetry of $\cK_0, \cK_1, \ldots, \cK_T$ to transform the absolute value of the sum into the sum of absolute values before decoupling the suprema.
Taking the defining expression of $\mu,\nu_1,\ldots,\nu_T$ into account finally leads to ${\rm gwce}(\Delta_a)$ being equal to the objective function of \eqref{MiniAbs}.
The announced result is consequently justified.
\epf

To make things less abstract,
we now introduce some particular model sets $\cK_0, \cK_1,\ldots,\cK_T$ based on approximation capabilities by linear subspaces $\cV_0,\cV_1,\ldots,\cV_T$ of $F$
with parameters $\eps_0,\eps_1,\ldots,\eps_T \hspace{-0.6mm}>\hspace{-0.6mm}0$.
As such,
let us consider
\be
\label{AppSet}
\cK_t = \{ g \in F: {\rm dist}(g,\cV_t) \le \eps_t \},
\qquad t = 0,1,\ldots,T.
\ee
For instance,
if $\cV_t = \{0\}$ for $t=1,\ldots,T$,
then the conditions $f_{t-1}-f_t \in \cK_t$
simply mean that $\|f_{t-1} - f_t\| \le \eps_t$,
i.e., that $f_{t-1}$ and $f_t$ are $\eps_t$-close.
Moreover,
the condition $f_0 \in \cK_0$ means that $f_0$ is $\eps_0$-close to being in the space $\cV_0$.
Such an assumption is made implicitly in many numerical procedures,
e.g. when devising quadrature formulas that are exact on a space of polynomials and whose accuracy depends on how close the function to integrate is from a polynomial.
Note that, as in the classical case where $\cK_1,\ldots,\cK_T$ are absent,
to avoid infinite worst-case errors,
it is assumed from the outset that 
$$
\cV_0 \cap \ker(\La_0) = \{0\},
\qquad \mbox{so that} \quad
n_0 := \dim(\cV_0) \le m_0.
$$
With the above approximability sets $\cK_0,\cK_1, \ldots, \cK_T$,
Theorem \ref{ThmLFAbs} becomes the following statement.

\bthm
\label{ThmLFApp}
Suppose that $\cK_0,\cK_1,\ldots,\cK_T$ are the approximability sets defined in \eqref{AppSet}.
If the quantity of interest $Q$ is a linear functional,
then a globally optimal recovery map is given by
$$
\Delta^{\rm opt}: y = [y_0; y_1; \ldots; y_T] \in \bR^{m} \mapsto 
\sum_{i=1}^{m_0} a^{\rm opt}_{0,i} y_{0,i}
+ \sum_{t=1}^T \sum_{i=1}^{m_t} a^{\rm opt}_{t,i} y_{t,i} \in \bR,
$$
where the vector $a^{\rm opt} = [a^{\rm opt}_0; a^{\rm opt}_1;\ldots; a^{\rm opt}_T] \in \bR^{m}$ is a solution to
\begin{align}
\minimize{a = [a_0; a_1;\ldots;a_T] \in \bR^{m} } 
& 
\left[ 
\left\| 
 Q - \sum_{t=0}^T \sum_{i=1}^{m_t} a_{t,i}\la_{t,i} \right\|_* \times \eps_0 
+ \sum_{s=1}^T
\left\| \sum_{t=s}^T \sum_{i=1}^{m_t} a_{t,i} 
\la_{t,i} \right\|_* \times \eps_s
\right]\\
\mbox{s.to} \quad & 
\left( \sum_{t=0}^T \sum_{i=1}^{m_t} a_{t,i} \la_{t,i} - Q \right)_{| \cV_0} = 0
\quad \mbox{and} \quad
\left( \sum_{t=s}^T \sum_{i=1}^{m_t} a_{t,i} \la_{t,i} \right)_{|\cV_s} =0
\mbox{ for all } s =1,\ldots,T.
\end{align}
\ethm

\bpf
The quantity to minimize in \eqref{MiniAbs} is made of several suprema,
each supremum taking the form $S = \sup_{g \in F} \{ |\eta(g)|: {\rm dist}(g,\cV) \le \eps \} $
for some linear functional $\eta$, linear subspace $\cV \inc F$, and parameter $\eps >0$.
Here, the constraint ${\rm dist}(g,\cV) \le \eps$ is equivalent to the existence of $v \in \cV$ such that $h:= g-v \in F$ satisfies $\|h\| \le \eps$.
It follows that the supremum $S$ can be written as
\begin{align*}
S & = \sup_{\substack{h \in F\\v \in \cV}} \{ |\eta(h+v)| :\|h\| \le \eps \}
= \sup_{h \in F} \{ |\eta(h)|: \|h\| \le \eps \}
+ \sup_{v \in \cV} \{ |\eta(v)| \} \\
& = \left\{
\bmx
+ \infty & \mbox{in case } \eta(v) \not= 0 \mbox{ for some }v \in \cV,\\
\|\eta\|_* \times \eps & \mbox{ in case } \mu(v) = 0 \mbox{ for all } v \in \cV. \quad \; \,
\emx
\right.
\end{align*}
Since we are interested in minimizing $S$,
we must place ourselves in the latter case and impose the constraint $\eta_{| \cV} = 0$.
Substituting the detailed expression of $\eta$ for each supremum yields the announced result. 
\epf

Still, it is not apparent that Theorem \ref{ThmLFApp} leads to an efficient construction of the vector $a^{\rm opt}$ because the dual norm on $F$ may not be fitted for computations.
For this to happen, we need to deal with spaces $F$ on a case-to-case basis.
For illustration, we shall take $F=C(\cX)$,
the space of continuous functions on a compact set $\cX$,
equipped with the norm $\|g\| = \max\{ |g(x)|, x \in \cX \}$.
As observation functionals,
we select point evaluations,
i.e., $\la_{t,i} = \delta_{x^{(t,i)}}$,
at points $x^{(t,i)}$ that are all distinct.
Thus, if $Q(f_0) = f_0(x)$ for some $x \not\in \{ x^{(t,i)}, t=0,1,\ldots, T, i=1,\ldots,m_t \}$
or $Q(f_0) = {\rm vol}(\cX)^{-1} \int_\cX f_0(x) dx$,
we can easily see (using Tietze extension theorem) that 
$$
\left\| 
 Q - \sum_{t=0}^T \sum_{i=1}^{m_t} a_{t,i}\la_{t,i} \right\|_* 
 =  1 + \sum_{t=0}^T \sum_{i=1}^{m_t} |a_{t,i}|,
 \qquad \quad
\left\| \sum_{t=s}^T \sum_{i=1}^{m_t} a_{t,i} 
\la_{t,i} \right\|_* 
= \sum_{t=s}^T \sum_{i=1}^{m_t} |a_{t,i}|,
\quad s=1,\ldots,T.
$$
With $(v_{t,1},\ldots,v_{t,n_t})$ denoting bases for $\cV_t$, $t=0,1,\ldots,T$,
we introduce the vector $b \in \bR^{n_0}$
and the matrices $M^{(s)} \in \bR^{n_s \times (m_s+\cdots+m_T)}$, $s=0,1,\ldots,T$,
with entries
$$
b_j =  v_{0,j}(x)
\qquad \mbox{and} \qquad
M^{(s)}_{j,(t,i)} = v_{s,j}(x^{(t,i)}).
$$
With this notation in place,
we can now state the following instantiation of Theorem \ref{ThmLFApp},
whose careful verification is left to the reader.

\bcor
\label{CorEstLFinCX}
Suppose that $\cK_0,\cK_1,\ldots,\cK_T \inc C(\cX)$ are the approximability sets defined in \eqref{AppSet}.
If the observation functionals $\la_{t,i}$ are evaluations at points $x^{(t,i)} \in \cX$ that are all distinct,
then a globally optimal recovery map for the evaluation of $f_0$ at new $x \in \cX$
is given by
$$
\Delta^{\rm opt}: y = [y_0; y_1; \ldots; y_T] \in \bR^{m} \mapsto 
\sum_{i=1}^{m_0} a^{\rm opt}_{0,i} y_{0,i}
+ \sum_{t=1}^T \sum_{i=1}^{m_t} a^{\rm opt}_{t,i} y_{t,i} \in \bR,
$$
where $a^{\rm opt} = [a^{\rm opt}_0; a^{\rm opt}_1; \ldots; a^{\rm opt}_T] \in \bR^{m}$ is a solution to
$$
\minimize{[a_0; a_1; \ldots; a_T] \in \bR^{m}} 
\sum_{t=0}^T (\eps_0+ \cdots+ \eps_t) \|a_t\|_1
\quad \mbox{s.to }
M^{(0)} [a_0; \ldots; a_T] = b,
\; M^{(s)} [a_s; \ldots; a_T ] = 0, \, s=1,\ldots,T.
$$
\ecor

It is worth mentioning that the above optimization program---a weighted $\ell_1$-minimization---can be recast as a standard-form linear program.
Furthermore, by freezing $a_1,\ldots,a_T$  to their optimal values,
we obtain a standard-form linear program with $n_0$ equality constraints in the variable $a_0 \in \bR^{m_0}$,
so one can output (e.g. by the simplex algorithm) a solution $a^{\rm opt}_0 \in \bR^{m_0}$ which is $n_0$-sparse.
Since the latter can be precomputed offline before the data $y \in \bR^m$ are collected,
it becomes unnecessary to collect the $y_{0,i}$'s for $i \not\in {\rm supp}(a^{\rm opt}_0)$.
This is potentially significant because the high-fidelity functionals $\la_{0,i}$ are typically expensive to evaluate.

\section{Globally Optimal Recovery in Hilbert Spaces}
\label{SecGlo}

This section deals with the second scenario,
where the quantity of interest is an arbitrary linear map $Q: F \to Z$---it is not restricted to be a functional anymore---but $F$ and $Z$ are restricted to be Hilbert spaces.
In close adherence to our multi-fidelity model,
there are only two model sets $\cK_0$ and $\cK_1$ here.
The latter are hyperellipsoids,
i.e., they have the form
\be
\label{2HyperEll}
\cK_0  = \{ g \in F: \|P_0 g \| \le \eps_0 \}
\qquad  \mbox{and} \qquad
\cK_1  = \{ g \in F: \|P_1 g \| \le \eps_1 \},
\ee
where the linear operators $P_0$ and $P_1$ also map into Hilbert spaces.
Note that the approximability sets introduced in \eqref{AppSet}
are covered in this setting,
as it suffices to choose $P_0 = P_{\cV_0^\perp}$ and $P_1 = P_{\cV_1^\perp}$,
i.e., the orthogonal projectors onto the orthogonal complements of $\cV_0$ and $\cV_1$, respectively.
For the full recovery problem,
i.e., when $Q = \Id_F$,
it is arguably natural to try and estimate the original $f_0$ by a constrained regularizer,
namely by the element $f^\tau_0$ obtained via
$$
[f_0^\tau; f_1^\tau]
:= \underset{ [f_0; f_1] \in F \times F}{\argmin \;}
\left[ (1-\tau) \|P_0 f_0\|^2 + \tau \|P_1(f_0 - f_1) \|^2 \right]
\qquad \mbox{s.to} \quad \La_0 f_0 = y_0, \; \La_1 f_1 = y_1
$$
for some parameter $\tau \in [0,1]$.
Rather than tuning this parameter by, say, cross-validation,
we shall uncover a principled way of selecting the optimal $\tau$.
By doing so, we will also guarantee that the associated constrained regularization map is globally optimal,
even among all possible recovery maps. 
The theorem stated below, 
which extends to the estimation of other quantities of interest besides $Q=\Id_F$,
follows from a recent result derived in~\cite[Section~2]{L1Err} on Optimal Recovery with a two-hyperellipsoid-intersection model set.
The key is to reduce our multi-fidelity setting to this related setting,
in the spirit of Section \ref{SecLF}.
Specifically, 
we reformulate the global worst-case error of any $\Delta: \bR^{m_0 + m_1} \to Z$ as
\be
\label{ReformGWCE_HGlo}
{\rm gwce}(\Delta) 
:= \sup_{\substack{\|P_0 f_0 \| \le \eps_0 \\ \|P_1(f_0 - f_1) \| \le \eps_1}} \|Q(f_0) - \Delta([\La_0 f_0; \La_1 f_1]) \|_Z
= \sup_{\substack{\|Rf\| \le 1\\ \|Sf\| \le 1 }} \|\wt{Q}(f) - \Delta( \La f) \|_Z,
\ee
where $f \in F \times F$ denotes again the compound vector
$$
f = \bbmx f_0 \\ \hline f_1 \ebmx,
$$
and where the Hilbert-valued linear maps $R$, $S$, $\La$, and $\wt{Q}$ are defined on $F \times F$ by
\be
\label{DefRS}
R(f) = \f{1}{\eps_0} P_0(f_0),
\qquad
S(f) = \f{1}{\eps_1} P_1(f_0 - f_1),
\qquad
\La f = \bbmx \La_0 f_0 \\ \hline \La_1 f_1 \ebmx,
\qquad
\wt{Q}(f) = Q(f_0). 
\ee
To avoid infinite worst-case errors, it is assumed from now on that 
$$
\ker(R) \cap \ker(S) \cap \ker(\La) = \{0\},
$$
or, when $\dim(F) = \infty$, that there is $\delta > 0$ such that $\max\{ \|Rh\|, \|Sh\| \} \ge \delta \|h\|$ for all $h \in \ker(\La)$.
In terms of our multi-fidelity framework,
this reads
\be
\label{CondPOP1La0La1}
[P_0(f_0) = 0, \; P_1(f_0-f_1) = 0, \; \La_0(f_0)=0, \; \La_1(f_1)=0  ]
\imp [f_0=0,\; f_1 =0].
\ee
In the cases $\ker(\La_0) \inc \ker(P_1)$ and $\ker(\La_0) \inc \ker(\La_1)$ featured in Section \ref{SecLoc}, 
it can be verified that \eqref{CondPOP1La0La1} holds if and only if both $\ker(P_0) \cap \ker(\La_0) = \{0\}$ and $\ker(P_1) \cap \ker(\La_1) = \{0\}$ hold.

\bthm
\label{ThmGlobOpt}
Suppose that $F$ is a Hilbert space and that the model sets $\cK_0,\cK_1 \inc F$
are the hyperellipsoidal sets defined in \eqref{2HyperEll}.
If the linear quantity of interest $Q: F \to Z$ maps into a Hilbert space, then a globally optimal recovery map is given by
$\Delta^{\rm opt}: y \in \bR^{m_0 + m_1}
\mapsto Q(f^{\tau^\sharp}_0) \in Z$,
where
$$
[f^{\tau^\sharp}_0; f^{\tau^\sharp}_1] = \underset{f = [f_0; f_1] \in F \times F}{\argmin \;}
\left[ (1-\tau^\sharp) \|P_0 f_0\|^2 + \tau^\sharp \|P_1(f_0 - f_1) \|^2 \right]
\qquad \mbox{s.to} \quad \La_0 f_0 = y_0, \; \La_1 f_1 = y_1. 
$$
Here, the regularization parameter is
$\tau^{\sharp} = c_1^{\sharp} / (c_0^{\sharp} + c_1^{\sharp})$, 
with $c_0^{\sharp}, c_1^{\sharp} \ge 0$ being solutions to
\begin{align}
\label{SDPGlo}
\minimize{c_0,c_1 \ge 0} &  \; [c_0 \eps_0^2 + c_1 \eps_1^2]\\
\nonumber
 \mbox{s.to } & \; c_0 \|P_0 f_0\|^2 + c_1 \|P_1(f_0 - f_1) \|^2 \ge \|Qf_0\|^2
\quad \mbox{for all } f_0 \in \ker(\La_0),
f_1 \in \ker(\La_1).
\end{align}
\ethm

\bpf
According to the reformulation \eqref{ReformGWCE_HGlo} of the global worst-case error of $\Delta: \bR^{m_0 + m_1} \to Z$
and invoking \cite[Theorem 1]{L1Err},
we know that an optimal recovery map is given by $\Delta^{\rm opt} = \wt{Q} \circ \Delta_{a^\sharp,b^\sharp}$, where
$$
\Delta_{a^\sharp,b^\sharp}:
y \in \bR^{m_0 + m_1}
\mapsto
\bigg(
 \underset{f \in F \times F}{\argmin \;}
a^\sharp \|Rf\|^2 + b^\sharp \|Sf\|^2   
\quad \mbox{s.to } \; \La f = y
\bigg) \in F \times F,
$$
and where $a^\sharp,b^\sharp \ge 0$
are solutions to 
$$
\minimize{a,b \ge 0}  \; (a  + b)
\quad \mbox{s.to }
a \|Rf \|^2 + b \|S f \|^2 \ge \|\wt{Q} f \|^2 
\mbox{ for all } f \in \ker \La.
$$
Making the change of optimization variables $c_0 = a/\eps_0^2$ and $c_1 = b/\eps_1^2$
while taking the expressions of $R$, $S$, $\La$, and $\wt{Q}$ into account,
we immediately see that the latter program reduces to \eqref{SDPGlo}.
Likewise,
we see that
$\Delta_{a^\sharp,b^\sharp}(y) =: [f^\sharp_0;f^\sharp_1]$ reduces to the minimizer  of 
$c_0^\sharp \|P_0 f_0\|^2 + c_1^\sharp \|P_1(f_0 - f_1) \|^2$---or equivalently of $(1-\tau^\sharp) \|P_0 f_0\|^2 + \tau^\sharp \|P_1(f_0 - f_1) \|^2$---subject to $\La_0 f_0 = y_0$ and $\La_1 f_1 = y_1$, i.e.,  
$[f^\sharp_0,f^\sharp_1] = [f^{\tau^\sharp}_0,f^{\tau^\sharp}_1]$,
so that $\Delta^{\rm opt}(y) = \wt{Q}([f^{\tau^\sharp}_0,f^{\tau^\sharp}_1]) = Q(f^{\tau^\sharp}_0)$, as announced.
\epf

\brk
Here are a few comments to put the theoretical results of Theorem \ref{ThmGlobOpt} into perspective.\vspace{-5mm}
\begin{enumerate}[a)]
\item Although not obvious at first sight,
the optimal recovery map $\Delta^{\rm opt} = \wt{Q} \circ \Delta_{a^\sharp,b^\sharp}$ is linear.
Indeed, the map $\Delta_{a^\sharp,b^\sharp}$ 
can be expressed (see \cite[Subsection 2.1]{L1Err}) as
$$ 
\Delta_{a^\sharp,b^\sharp} = \La^\dagger - \big[ a^\sharp R_{\cN}^* R_{\cN} + b^\sharp S_{\cN}^* S_{\cN} \big]^{-1} \big( a^\sharp R_{\cN}^* R + b^\sharp S_{\cN}^* S \big) \La^\dagger,
$$
where $R_{\cN}$ and $S_{\cN}$ denote the restrictions of $R$ and $S$ to $\cN = \ker(\La)$.

\item Theorem \ref{ThmGlobOpt} is valid for $\dim(F)=\infty$ (since the underpinning result from \cite{L1Err}~is),
but it is in finite dimensions that 
\eqref{SDPGlo} becomes a manageable semidefinite program.
Indeed, if $(h^0_1,\ldots,h^0_{k_0})$ and $(h^1_1,\ldots,h^1_{k_1})$ are bases for $\ker(\La_0)$ and $\ker(\La_1)$ and if $H_0: \bR^{k_0} \to \ker(\La_0)$ and $H_1: \bR^{k_1} \to \ker(\La_1)$ represent the transformations defined by $H_0(u) = \sum_{i=1}^{k_0} u_i h^0_i$ and
$H_1(v) = \sum_{i=1}^{k_1} v_i h^1_i$,
then we can reformulate the constraint 
$c_0 \|P_0 f_0\|^2 + c_1 \|P_1(f_0 - f_1) \|^2 \ge \|Qf_0\|^2$
for all $f_0 \in \ker(\La_0)$ and 
$f_1 \in \ker(\La_1)$
as
$$
\bbmx
H_0^*(c_0 P_0^* P_0 + c_1 P_1^* P_1 - Q^* Q) H_0
& \vline &
H_0^*(c_1 P_1^* P_1) H_1\\
\hline
H_1^*(c_1 P_1^* P_1) H_0
& \vline & 
H_1^* (c_1 P_1^* P_1) H_1
\ebmx \succeq 0.
$$
Note that the square matrix above has size $k_0+k_1 \ge 2 \dim(F) - m_0 - m_1$,
independently of the dimension of $Z$.
Intuitively, if $Z$ is low dimensional,
e.g. $Z = \bR$ when $Q$ is a linear functional,
then we should be able to lower the computational cost of producing an optimal recovery map,
but this does not follow from Theorem~\ref{ThmGlobOpt}.

\item  Theorem \ref{ThmGlobOpt} substantially relies on a result from \cite{L1Err},
itself relying on the so-called S-procedure,
and in particular on its `no-linear-term' version.
Loosely speaking,
given quadratic functions $q_0,q_1, \ldots, q_K$~on~$\bK^n$,
the S-procedure---which is surveyed in \cite{PolTer}---relates the condition
[$q_0(x) \le 0$ whenever $q_1(x) \le~0, \ldots, q_K(x) \le 0$]
to the condition of existence of $c_1,\ldots,c_K \ge 0$ such that $q_0 \le c_1 q_1 + \cdots + c_K q_K$.
Clearly, the latter condition always implies the former condition.
The S-procedure is said to be exact if, conversely, the former implies the latter.
Exactness occurs, under some strict feasibility conditions,
when $K=1$ (this is Yakubovich's S-Lemma) 
and when $K=2$ if $\bK = \bC$ (see \cite{BecEld2}).
It also occurs when $K=2$ if $\bK = \bR$
provided the quadratic functions are all of the form $q_k(x) = \langle A_k x, x \rangle + \alpha_k$,
i.e., they do not feature any linear term $\langle a_k,x \rangle$.
Established by Polyak in \cite{Pol},
this is what we called the `no-linear-term' S-procedure.
Not only is this version generally invalid for $K>2$,
but the result from~\cite{L1Err} is invalid, too,
when more than two hyperellipsoidal model sets are involved.
Consequently, unless a more restrictive setting is enforced,
we cannot extend Theorem \ref{ThmGlobOpt} to the situation considered in Section \ref{SecLF},
 where we dealt with model sets $\cK_0,\cK_1,\ldots, \cK_T$ with $T>1$.

\end{enumerate}
\erk

\section{Chebyshev Center for the Intersection of Two Hyperellipsoids}
\label{SecCCtwoHyp}

Still within Hilbert spaces,
our next goal is to tackle the more subtle problem of local optimality.
As was just made apparent in the global setting,
our multi-fidelity model is tied to the two-hyperellipsoid-intersection model,
so we start with some investigations about this general model, bearing in mind that the local setting was not considered in \cite{L1Err}.
Readers mainly interested  in the multi-fidelity aspect can skip to Section \ref{SecLoc} and bypass the technicalities of the present section.
Here, we deal with Hilbert-valued linear maps $R$, $S$, $\La$, and $\wt{Q} $ which are defined on a common Hilbert space~$F$
and which do not have to be taken as in \eqref{DefRS}.
Fixing $y \in {\rm ran}(\La)$,
we aim at minimizing over all~$z$ the local worst-case error 
$$
{\rm lwce}_y(z) = \sup_{\substack{\|Rf\| \le 1, \|Sf\| \le 1\\ \La f = y}} \|\wt{Q}(f) - z\|.
$$
In other words, we are looking for a Chebyshev center of the linear image of a sliced intersection of~two hyperellispoids,
namely for a center $z^{\rm opt}$ of a ball with smallest radius that contains the set $\wt{Q}(\{ f \in F: \|Rf\| \le 1, \|Sf\| \le 1, \La f = y \})$.
In this situation,
classical results guarantee the existence and~uniqueness of such a Chebyshev center.
However,
determining the Chebyshev center/radius is an arduous task in general.
In the closely related problem of Optimal Recovery from inaccurate data,
a  semidefinite-programming-based upper bound for the Chebyshev radius was proposed in~\cite{BecEld} when $Q = \Id_F$, 
along with a candidate Chebyshev center.
The genuine Chebyshev center was obtained more directly in \cite{FouLia} but under a specific assumption and still when $Q = \Id_F$.
It was later observed that the aforementioned  candidate and genuine Chebyshev centers actually coincide,
see~\cite{FouLia-CAMDA}.
There, the upper bound was also generalized to linear quantities of interest $Q \not= \Id_F$.
All these findings can now be retrieved from the theorem below.
It consists of three parts:
the first part provides a computable upper bound for the Chebyshev radius (incidentally involving constrained regularizers);
the second part uncovers an orthogonality condition under which this upper bound agrees with the true Chebyshev radius;
the third part shows conversely that the upper bound and the true Chebyshev radius can only be equal under this orthogonality condition.
The precise statement of the theorem involves
the following notation to be used throughout this section: given $y \in {\rm ran}(\La)$ and $\tau \in [0,1]$, \vspace{-5mm}
\begin{itemize}
\item $f^\tau \in F$ is the constrained minimizer defined as
$$
f^\tau := \underset{{f \in F}}{\argmin} \; \big[ (1-\tau) \|Rf\|^2 + \tau \|S f\|^2 \big]
\qquad \mbox{s.to} \quad \La f = y;
$$ 
\item $\la(\tau)$ denotes the largest eigenvalue of 
$[(1-\tau) R_\cN^* R_\cN + \tau S_\cN^* S_\cN]^{-1} \wt{Q}_\cN^* \wt{Q}_\cN$,  $\cN := \ker(\La)$, i.e.,
$$
\la(\tau) := \la_{\max}([(1-\tau) R_\cN^* R_\cN + \tau S_\cN^* S_\cN]^{-1} \wt{Q}_\cN^* \wt{Q}_\cN);
$$
\item $h^\tau \in \cN$ is defined implicitely via an eigenequation
and a normalization, namely
\begin{align*}
\wt{Q}_\cN^* \wt{Q}_\cN h^\tau 
& = \la(\tau) [(1-\tau) R_\cN^* R_\cN + \tau S_\cN^* S_\cN] h^\tau,\\ 
\|\wt{Q} h^\tau\|^2 
& = \la(\tau) \times  \big( 1 - ( (1-\tau) \|R f^\tau\|^2 + \tau \| S f^\tau \|^2 ) \big).
\end{align*}
\end{itemize}
The three-part theorem, whose upcoming proof is technically involved, reads as follows.

\bthm
\label{ThmCCtwoHyp} 
Let $R,S,\wt{Q}$ be Hilbert-valued linear maps defined on a common Hilbert space $F$.
Let also $\La: F \to \bR^m$ be another linear map and let $y \in \bR^m$.\\
(i) The squared Chebyshev radius $\rho^2$ of the set $\wt{Q}(\{ f \in F: \|R f \| \le 1, \|S f\| \le 1, \La f = y \})$ is upper-bounded as
$$
\rho^2 
\le \inf_{\tau \in [0,1]} \sup_{\substack{\|Rf\| \le 1, \|Sf\| \le 1\\ \La f = y}} \|\wt{Q}(f) - \wt{Q}(f^\tau)\|^2
\le \inf_{\tau \in [0,1]} \| \wt{Q} h^\tau \|^2.
$$
(ii) For the minimizer $\tau^y$ of $\|Q h^\tau\|^2$ over $\tau \in [0,1]$, if the orthogonality condition
\be
\label{OrthCond}
\langle R f^{\tau^y}, R h^{\tau^y} \rangle = 0
\qquad \mbox{and/or} \qquad 
\langle S f^{\tau^y}, S h^{\tau^y} \rangle = 0
\ee
holds,
then the squared Chebyshev radius agrees with the upper bound, i.e.,
\be
\label{SqChe=}
\rho^2 = \| \wt{Q} h^{\tau^y} \|^2
= \inf_{\tau \in [0,1]} \| \wt{Q} h^\tau \|^2.
\ee
(iii) Conversely, if the equality \eqref{SqChe=} between squared Chebyshev radius and upper bound occurs,
then the orthogonality condition \eqref{OrthCond} holds.
\ethm

Each part of the theorem will be proved separately.
For the first part, we make use of the following lemma,
which provides a fresh look at the one-hyperellipsoid model.
The result was obtained in \cite{BCDDPW}
when $T = P_{\cV^\perp}$ is the orthogonal projector onto the orthogonal complement of a linear space $\cV \inc F$ and when $Q = \Id_F$,
in which case the largest eigenvalue in the statement below reduces to $1/\sigma_{\min}(T_\cN)^2$.
When $T = P_{\cV^\perp}$ but $Q \not= \Id_F$,
it can also be found as \cite[Theorem 10.2]{BookDS}, albeit with the largest eigenvalue expressed in a different manner.

\blem
\label{LemCC1Hyp}
Let $T,\wt{Q}$ be Hilbert-valued linear maps defined on a common Hilbert space $F$.
Let also $\La: F \to \bR^m$ be another linear map whose null space $\cN = \ker(\La)$ satistfies $\ker(T) \cap \cN = \{0\}$.
Given $y \in \bR^m$,
the Chebyshev center of the set $\wt{Q}(\{ f \in F: \|T f \| \le 1, \La f = y \})$ is $\wt{Q}(\wh{f})$,
where
$$
\wh{f} := \underset{f \in F}{\argmin } \; \|Tf\|^2 \qquad \mbox{s.to} \quad \La f = y
$$
and its squared Chebyshev radius is equal to 
\be
\label{Normahhat}
\|\wt{Q} \wh{h}\|^2 = 
\la_{\max}\big( [T_\cN^* T_\cN ]^{-1} \wt{Q}_\cN^* \wt{Q}_\cN \big) \times \big( 1 - \|T \wh{f} \|^2 \big),
\ee
where $\wh{h} \in \cN$ is a leading eigenvector of $[T_\cN^* T_\cN ]^{-1} \wt{Q}_\cN^* \wt{Q}_\cN $ normalized so that $\|T \wh{h} \|^2 = 1 - \| T \wh{f} \|^2$.
\elem

\bpf
Let us first quickly justify the equality in \eqref{Normahhat}, which relies on the properties of $\wh{h}$ via
\begin{align*}
\|\wt{Q} \wh{h}\|^2  
& = \langle \wt{Q}_\cN^* \wt{Q}_\cN \wh{h}, \wh{h} \rangle
= \langle  [T_\cN^* T_\cN ]^{-1} \wt{Q}_\cN^* \wt{Q}_\cN \wh{h}, [T_\cN^* T_\cN ] \wh{h} \rangle
= \langle  \la_{\max}([T_\cN^* T_\cN ]^{-1} \wt{Q}_\cN^* \wt{Q}_\cN ) \, \wh{h}, [T_\cN^* T_\cN ] \wh{h} \rangle\\
& = \la_{\max}([T_\cN^* T_\cN ]^{-1} \wt{Q}_\cN^* \wt{Q}_\cN ) \times \|T \wh{h}\|^2
 = \la_{\max}([T_\cN^* T_\cN ]^{-1} \wt{Q}_\cN^* \wt{Q}_\cN ) \times \big( 1 - \| T \wh{f} \|^2 \big).
\end{align*}
The equality in \eqref{Normahhat} being put aside, our remaining objective is to prove that, for any $z \in {\rm ran}(\wt{Q})$,
\be
\label{Obj}
\sup_{\substack{\|Tf \|^2 \le 1\\ \La f = y}} \|\wt{Q}f-z\|^2
\ge   
\sup_{\substack{\|Tf \|^2 \le 1\\ \La f = y}} \|\wt{Q} f- \wt{Q} \wh{f}\|^2 
=  \|\wt{Q} \wh{h}\|^2.
\ee
We start with the equality on the right of \eqref{Obj}.
To this end, 
it is crucial to observe that $\wh{f}$ is characterized by  $\La \wh{f} = y$ and $\langle T \wh{f}, T h \rangle = 0$ for all $h \in \cN$,
the latter property stemming from $\|T \wh{f}\|^2 \le \|T (\wh{f}+h)\|^2$ for all $h \in \cN$.
Thus, reparametrizing any 	$f \in F$ such that $\La f =y$ as $f = \wh{f} + h$ with $h \in \cN$,
we see that the constraint $\|Tf \|^2 \le 1$,
i.e., $\|T(\wh{f} + h) \|^2 \le 1$,
is equivalent to $\|T\wh{f}\|^2 + \| T h \|^2 \le 1$.
Therefore, 
$$
\sup_{\substack{\|Tf \|^2 \le 1\\ \La f = y}} \|\wt{Q}f- \wt{Q} \wh{f}\|^2  
= \sup_{\substack{\|Th \|^2 \le 1 - \|T \wh{f} \|^2 \\ h \in \cN}} \|\wt{Q} h\|^2.
$$
At this point, we claim that the latter supremum is equal to $\|\wt{Q} \wh{h}\|^2$.
Indeed,
for any $h \in \cN$ satisfying $\|Th\|^2 \le 1 - \|T \wh{f}\|^2$, 
by introducing $u := [T_\cN^* T_\cN ]^{1/2} h \in \cN$ and using the shorthand $\la_{\max} := \la_{\max}\big( [T_\cN^* T_\cN ]^{-1} \wt{Q}_\cN^* \wt{Q}_\cN \big) = 
\la_{\max} \big( [T_\cN^* T_\cN ]^{-1/2} \wt{Q}_\cN^* \wt{Q}_\cN [T_\cN^* T_\cN ]^{-1/2} \big)$, 
we remark that
\begin{align*}
\|\wt{Q} h\|^2 & = \langle \wt{Q}_\cN^* \wt{Q}_\cN h, h \rangle
= \langle [T_\cN^* T_\cN ]^{-1/2} \wt{Q}_\cN^* \wt{Q}_\cN [T_\cN^* T_\cN ]^{-1/2} u , u \rangle \\
& \le \la_{\max}  \times \|u\|^2
= \la_{\max}  \times \|Th \|^2
\le \la_{\max}  \times \big( 1 - \|T \wh{f} \|^2 \big),
\end{align*}
with equality occurring throughout if and only if
$[T_\cN^* T_\cN ]^{-1/2} \wt{Q}_\cN^* \wt{Q}_\cN [T_\cN^* T_\cN ]^{-1/2} u = \la_{\max} \, u$,
i.e., $\wt{Q}_\cN^* \wt{Q}_\cN h = \la_{\max} \, [T_\cN^* T_\cN ] \, h$,
and $\|T h \|^2 = 1 - \|T \wh{f} \|^2$---
in short, if and only if $h$ equals (one of the)~$\wh{h}$.

Continuing with the inequality on the left of \eqref{Obj}, let us consider any $z \in {\rm ran}(\wt{Q})$.
By noticing that $\|T(\wh{f} \pm \wh{h})\|^2 = 1$ and that $\La (\wh{f} \pm \wh{h}) = y$,
we obtain
\begin{align*}
\sup_{\substack{\|Tf \|^2 \le 1\\ \La f = y}} \|\wt{Q} f-z\|^2
& \ge \max_\pm \|\wt{Q} (\wh{f} \pm \wh{h}) - z\|^2
\ge \f{1}{2}  \|\wt{Q} \wh{f} - z + \wt{Q} \wh{h}\|^2
+ \f{1}{2}  \|Q\wh{f} - z - \wt{Q}\wh{h}\|^2\\
& = \|\wt{Q}\wh{f} - z\|^2 +  \|\wt{Q}\wh{h}\|^2
\ge \|\wt{Q}\wh{h}\|^2,
\end{align*}
which is the desired inequality.
The proof is now complete.
\epf

\bpf[Proof of Theorem \ref{ThmCCtwoHyp}, Part (i)]
For a fixed $\tau \in [0,1]$,
the squared Chebyshev radius $\rho^2$, 
i.e., the squared minimal local worst-case error,
is evidently upper-bounded by the squared local worst-case error evaluated at $\wt{Q}(f^\tau)$.
We bound the latter as 
$$
\sup_{\substack{\|Rf\|^2 \le 1, \|Sf\|^2 \le 1\\ \La f = y}} \|\wt{Q}(f) - \wt{Q}(f^\tau)\|^2
\le 
\sup_{\substack{(1-\tau)\|Rf\|^2 + \tau \|Sf\|^2 \le 1\\ \La f = y}} \|\wt{Q}(f) - \wt{Q}(f^\tau)\|^2.
$$
If we set $T := [(1-\tau) R^* R + \tau S^* S]^{1/2}$,
so that $\|Tf\|^2 = (1-\tau) \|Rf\|^2 + \tau \|Sf\|^2$
for all $f \in F$,
we can invoke Lemma \ref{LemCC1Hyp}---and its proof, see \eqref{Obj} in particular---to state that the latter supremum is equal to
$\|\wt{Q} \wh{h}\|^2$,
where $\wh{h}$ is easily identified in our context with $h^\tau$.
We have therefore obtained 
$$
\rho^2 \le
\sup_{\substack{\|Rf\|^2 \le 1, \|Sf\|^2 \le 1\\ \La f = y}} \|\wt{Q}(f) - \wt{Q}(f^\tau)\|^2
\le \|\wt{Q} h^\tau \|^2.
$$
The required result simply follows by taking the infimum over $\tau \in [0,1]$.
\epf

We now turn to the second part of the theorem,
for which it is worth isolating another lemma.
After uncovering the derivatives of some $\tau$-dependent quantities,
the lemma brings forward a crucial piece of information about the minimizer $\tau^y$ of $\|Q h^\tau\|^2$.

\blem
\label{LemDiffGH}
The derivatives of the functions $G$ and $H$ defined below are
\begin{align*}
\mbox{if } G(\tau) & := (1-\tau) \|R f^\tau\|^2 + \tau \|S f^\tau\|^2,
& \mbox{then} \qquad G'(\tau) &= \|Sf^\tau\|^2 - \|R f^\tau\|^2,\\
\mbox{if } H(\tau) & := 1/\la(\tau),
& \mbox{then} \qquad H'(\tau) & = \f{\|Sh^\tau\|^2 - \|R h^\tau\|^2}{\|\wt{Q} h^\tau \|^2}.
\end{align*}
As a consequence, the minimizer $\tau^y \in (0,1)$ of $\|Q h^\tau\|^2$  satisfies
$$
\|R f^{\tau^y} \|^2 + \|R h^{\tau^y} \|^2
= \|S f^{\tau^y} \|^2 + \|S h^{\tau^y} \|^2
=1.
$$
\elem

\bpf
For $\sigma,\tau \in [0,1]$, 
the minimality property of $f^\sigma$ yields
\begin{align*}
 (1-\sigma) \|R f^\sigma\|^2 + \sigma \|S f^\sigma \|^2
& \le (1-\sigma) \|R f^\tau\|^2 + \sigma \|S f^\tau \|^2\\
& = (1-\tau) \|R f^\tau\|^2 + \tau \|S f^\tau \|^2
+ (\sigma-\tau) (\|S f^\tau\|^2 - \|R f^\tau\|^2).
\end{align*}
In terms of the function $G$, the latter reads
$G(\sigma) \le G(\tau) + (\sigma-\tau) (\|S f^\tau\|^2 - \|R f^\tau\|^2)$.
Likewise, by exchanging the roles of $\tau$ and $\sigma$, we also have
$G(\tau) \le G(\sigma) + (\tau-\sigma) (\|S f^\sigma\|^2 - \|R f^\sigma\|^2)$.
We combine these two inequalities in the form
$$
(\sigma-\tau) (\|S f^\sigma\|^2 - \|R f^\sigma\|^2)
\le G(\sigma) - G(\tau) 
\le (\sigma-\tau) (\|S f^\tau\|^2 - \|R f^\tau\|^2).
$$
This immediately implies that the function $G$ is differentiable with
$$
G'(\tau) = \|S f^\tau\|^2 - \|R f^\tau\|^2.
$$
Next, the eigenequation defining $h^\tau$ can be written as
$ [(1-\tau) R_\cN^* R_\cN + \tau S_\cN^* S_\cN] h^\tau = H(\tau) \wt{Q}_\cN^* \wt{Q}_\cN h^\tau$,
and differentiating it with respect to $\tau$, we obtain
$$
(- R_\cN^* R_\cN +  S_\cN^* S_\cN) h^\tau + [(1-\tau) R_\cN^* R_\cN + \tau S_\cN^* S_\cN]\f{d h^\tau}{d \tau}
= H'(\tau) \wt{Q}_\cN^* \wt{Q}_\cN h^\tau + H(\tau) \wt{Q}_\cN^* \wt{Q}_\cN \f{d h^\tau}{d \tau}.
$$
Taking the inner product with $h^\tau$, while noticing that 
\begin{align*}
\left\langle [(1-\tau) R_\cN^* R_\cN + \tau S_\cN^* S_\cN]\f{d h^\tau}{d \tau}, h^\tau \right\rangle
& = \left\langle \f{d h^\tau}{d \tau}, [(1-\tau) R_\cN^* R_\cN + \tau S_\cN^* S_\cN] h^\tau \right\rangle \\
& = \left\langle \f{d h^\tau}{d \tau}, H(\tau) \wt{Q}_\cN^* \wt{Q}_\cN h^\tau \right\rangle 
= H(\tau) \left\langle \wt{Q}_\cN^* \wt{Q}_\cN  \f{d h^\tau}{d \tau},  h^\tau \right\rangle,
\end{align*}
we arrive at the identity $\langle (-R_\cN^* R_\cN + S_\cN^* S_\cN) h^\tau, h^\tau \rangle = H'(\tau) \langle \wt{Q}_\cN^* \wt{Q}_\cN h^\tau, h^\tau \rangle$,
which also reads $-\|R h^\tau\|^2 + \|S h^\tau\|^2 = H'(\tau) \|\wt{Q} h^\tau \|^2$.
This yields the announced expression for $H'(\tau)$.

Finally, from the defining property of $h^\tau$,
we point out that $\| \wt{Q} h^\tau \|^2 = (1-G(\tau))/H(\tau)$
and remark that the derivative of $H$ can also be expressed as $H'(\tau) = (\|S h^\tau\|^2 - \|R h^\tau\|^2 ) H(\tau)/(1-G(\tau))$.
It follows that,
at the minimizer $\tau^y \in (0,1)$ of $\|\wt{Q} h^\tau \|$, we have 
\begin{align}
\label{Der=0}
0 & = -G'(\tau^y) H(\tau^y) - (1-G(\tau^y)) H'(\tau^y)\\
\nonumber 
& = \big( \|R f^{\tau^y} \|^2 - \|S f^{\tau^y}\|^2
+ \|R h^{\tau^y}\|^2 - \|Sh^{\tau^y}\|^2 \big)
\times H(\tau^y).
\end{align}
Simplifying by $H(\tau^y) \not= 0$ and rearranging yields $\|R f^{\tau^y} \|^2 + \|R h^{\tau^y}\|^2 = \|S f^{\tau^y}\|^2 + \|Sh^{\tau^y}\|^2$.
With~$\kappa$ denoting this common value,
we see that
\begin{align*}
\kappa &= 
(1-\tau^y) \big(\|R f^{\tau^y}\|^2 + \|R h^{\tau^y}\|^2 \big) + \tau^y \big( \|S f^{\tau^y}\|^2 + \|S h^{\tau^y}\|^2 \big) \\
& = \big( (1-\tau^y)\|R f^{\tau^y}\|^2 + \tau^y \|S f^{\tau^y}\|^2 \big) + \big( (1-\tau^y)\|R h^{\tau^y}\|^2 + \tau^y \|S h^{\tau^y}\|^2 \big).
\end{align*}
Taking into account that
\begin{align*}
(1-\tau^y)\|R h^{\tau^y}\|^2 + \tau^y \|S h^{\tau^y}\|^2
& = \langle [(1-\tau^y) R_\cN^* R_\cN + \tau S_\cN^* S_\cN] h^{\tau^y}, h^{\tau^y} \rangle
= 
\langle (1/\la(\tau^y)) \wt{Q}_\cN^* \wt{Q}_\cN h^{\tau^y} , h^{\tau^y} \rangle\\
& = (1/\la(\tau^y)) \|\wt{Q} h^{\tau^y}\|^2  
= 1 - ( (1-\tau^y) \|R f^{\tau^y}\|^2 + \tau^y \| S f^{\tau^y} \|^2 ) ,
\end{align*}
we derive that $\kappa = 1$.
The desired result is now proved.
\epf

\bpf[Proof of Theorem \ref{ThmCCtwoHyp}, Part (ii)]
Together with $\|R f^{\tau^y} \|^2 + \|R h^{\tau^y} \|^2
= \|S f^{\tau^y} \|^2 + \|S h^{\tau^y} \|^2
=1$,
the additional orthogonality conditions $\langle R f^{\tau^y}, R h^{\tau^y} \rangle  = 0$
and $\langle S f^{\tau^y}, S h^{\tau^y} \rangle  = 0$---note that one follows from the other according to the characterization $((1-\tau^y)R^* R + \tau S^* S) f^{\tau^y} \perp \ker \La$---imply
$$
\|R(f^{\tau^y} \pm h^{\tau^y}) \|^2 
= 1
\qquad \mbox{and} \qquad
\|S(f^{\tau^y} \pm h^{\tau^y}) \|^2 
= 1.
$$
From here, it follows that, for any $z \in {\rm ran}(\wt{Q})$,
\begin{align*}
{\rm lwce}_y(z)^2
& = \sup_{\substack{\|Rf\| \le 1, \|Sf\| \le 1\\ \La f = y}} \|\wt{Q}(f) - z\|^2
\ge \max_{\pm} \|\wt{Q}(f^{\tau^y} \pm h^{\tau^y}) - z\|^2\\
& \ge \f{1}{2} \| \wt{Q} f^{\tau^y} - z + \wt{Q} h^{\tau^y} \|^2
+  \f{1}{2} \| \wt{Q} f^{\tau^y} - z - \wt{Q} h^{\tau^y} \|^2
=  \| \wt{Q} f^{\tau^y} - z\|^2
+  \| \wt{Q} h^{\tau^y}\|^2\\
& \ge \| \wt{Q} h^{\tau^y} \|^2.
\end{align*}
This lower bound is the same as the upper bound from Part (i),
so that Part (ii) is now proved.
\epf

\bpf[Proof of Theorem \ref{ThmCCtwoHyp}, Part (iii)]
Let us now assume that the squared Chebyshev radius coincides with the upper bound from Part (i),
i.e., that the minimum squared local worst-case error equals $\|\wt{Q} h^{\tau^y}\|^2$.
In this case, we have
$$
\|\wt{Q} h^{\tau^y}\|^2 
\le {\rm lwce}_y(\wt{Q} f^{\tau^y})
= \sup_{\substack{\|Rf\| \le 1, \|Sf\| \le 1 \\ \La f = y}} \|\wt{Q}f - \wt{Q}f^{\tau^y} \|^2
= \sup_{\substack{\|R(f^{\tau^y} + h)\| \le 1 \\ 
\|S(f^{\tau^y} + h)\| \le 1 \\ h \in \cN}} \|\wt{Q} h \|^2.
$$
If the latter supremum is attained at some $\ol{h} \in \cN$, we have $ \|\wt{Q} h^{\tau^y}\|^2 \le \|\wt{Q} \ol{h}\|^2$.
Using the shorthand $T := [(1-\tau^y) R^* R + \tau^y S^* S]^{1/2} $, the following chain of inequalities makes use of $T^* T f^{\tau^y} \perp \cN$, of $\|R(f^{\tau^y} + \ol{h})\| \le 1$,  $\|S(f^{\tau^y} + \ol{h})\| \le 1$, 
and of the normalization of $h^{\tau^y}$:
\begin{align*}
\| T \ol{h} \|^2
& = \langle T^* T \ol{h}, \ol{h} \rangle
= \langle T^* T(f^{\tau^y} + \ol{h}), (f^{\tau^y} + \ol{h}) \rangle - \langle T^* T f^{\tau^y}, f^{\tau^y} \rangle\\ 
& = \langle [(1-\tau^y) R^* R  + \tau^y S^* S] (f^{\tau^y} + \ol{h}), (f^{\tau^y} + \ol{h}) \rangle
- \langle [(1-\tau^y) R^* R  + \tau^y S^* S] f^{\tau^y}, f^{\tau^y} \rangle\\
& = \big( (1-\tau^y) \|R(f^{\tau^y} + \ol{h}) \|^2 + \tau^y \| S(f^{\tau^y} + \ol{h}) \|^2 \big) 
- \big( (1-\tau^y) \|R f^{\tau^y} \|^2 + \tau^y \| S f^{\tau^y}\|^2 \big)\\
& \le \big( (1-\tau^y) \times 1 + \tau^y \times 1 \big)
- \big( (1-\tau^y) \|R f^{\tau^y} \|^2 + \tau^y \| S f^{\tau^y} \|^2 \big)\\
& = 1 - \big( (1-\tau^y) \|R f^{\tau^y}  \|^2 + \tau^y \| S f^{\tau^y} \|^2 \big) 
= \f{1}{\la(\tau^y)} \|\wt{Q}h^{\tau^y} \|^2\\
& \le \f{1}{\la(\tau^y)} \|\wt{Q} \ol{h} \|^2
= \f{1}{\la(\tau^y)} \langle
[T_\cN^* T_\cN]^{-1/2} \wt{Q}_\cN^* \wt{Q}_\cN [T_\cN^* T_\cN]^{-1/2} [T_\cN^* T_\cN]^{1/2} \ol{h}, 
[T_\cN^* T_\cN]^{1/2} \ol{h} \rangle \\
& \le \f{1}{\la(\tau^y)} \la_{\max} \big( [T_\cN^* T_\cN]^{-1/2} \wt{Q}_\cN^* \wt{Q}_\cN [T_\cN^* T_\cN]^{-1/2} \big)
\| [T_\cN^* T_\cN]^{1/2} \ol{h} \|^2 \\
& = \| T_\cN \ol{h} \|^2.
\end{align*}
Since the leftmost term and the rightmost term are the same, equality must hold all the way through,
implying that $\|R(f^{\tau^y} + \ol{h}) \|^2 = 1$, $\|S(f^{\tau^y} + \ol{h}) \|^2 = 1$,
and that $[T_\cN^* T_\cN]^{1/2} \ol{h}$ is a leading eigenvalue of $[T_\cN^* T_\cN]^{-1/2} \wt{Q}_\cN^* \wt{Q}_\cN [T_\cN^* T_\cN]^{-1/2}$,
hence $\ol{h} = h^\tau$.
We thus arrived at $\|R(f^{\tau^y} + h^{\tau^y}) \|^2 = 1$.
Recalling from Lemma \ref{LemDiffGH} that $\|Rf^{\tau^y} \|^2  + \|R h^{\tau^y} \|^2 = 1$, 
we conclude that $\langle R f^{\tau^y}, R h^{\tau^y} \rangle =0$, as desired.
\epf

\brk
Although Theorem \ref{ThmCCtwoHyp} may appear like an abstract statement distant from practical considerations,
our candidate locally optimal recovery map,
as well as a locally optimal one,
can actually be efficiently computed,
as revealed in the comments below.\vspace{-5mm}
\begin{enumerate}[a)]

\item
The upper bound $\|Qh^{\tau^y}\|^2$ from Theorem \ref{ThmCCtwoHyp}, Part (i),
is the optimal value of the semidefinite program
\begin{align}
\label{SDP2Hyper}
\minimize{b, c, d \ge 0} & \; 
c (1 - \|R \La^\dagger y \|^2) + d( 1- \|S \La^\dagger y \|^2 ) + b\\
\nonumber
\mbox{s.to } &  \;  c  R_\cN^* R_\cN + d S_\cN^* S_\cN  \succeq \wt{Q}_\cN^* \wt{Q}_\cN\\
\nonumber
\mbox{and } &  \;
\bbmx
c R_\cN^* R_\cN + d S_\cN^* S_\cN  & \vline & 
c R_\cN^* R \La^\dagger y + d S_\cN^* S \La^\dagger y \\
\hline
\big( c R_\cN^* R \La^\dagger y + d S_\cN^* S \La^\dagger y \big)^* & \vline & b
\ebmx \succeq 0. 
\end{align}
Moreover, if $b^y,c^y,d^y \ge 0$
are minimizers, then the optimal parameter $\tau^y$ is given by $d^y/(c^y + d^y)$.

To justify this claim, we first look at the one-hyperellipsoid situation. 
By Lemma \ref{LemCC1Hyp},
the squared Chebyshev radius is
$\rho^2 = \la_{\max} \, \big( 1-\|T \wh{f}\|^2 \big)$, 
where $\la_{\max} := \la_{\max}([T_\cN^* T_\cN]^{-1/2} \wt{Q}_\cN^* \wt{Q}_\cN [T_\cN^* T_\cN]^{-1/2})$
can be expressed as
$$
\la_{\max} = \inf_{a \ge 0} \big\{  a: \;  a \Id_\cN \succeq [T_\cN^* T_\cN]^{-1/2} \wt{Q}_\cN^* \wt{Q}_\cN [T_\cN^* T_\cN]^{-1/2} \big\}
= \inf_{a \ge 0} \big\{  a: \;  a T_\cN^* T_\cN \succeq  \wt{Q}_\cN^* \wt{Q}_\cN  \big\}.
$$
Keeping in mind that $T \wh{f} = y$ and that $T^* T \wh{f} \perp \cN$---which yields $\wh{f} = \La^\dagger y - [T_\cN^* T_\cN]^{-1} T_\cN^* T \La^\dagger y$---we note that
\begin{align*}
\|T \wh{f}&  \|^2 
= \langle T \wh{f}, T \La^\dagger y + T( \wh{f} -\La^\dagger y) \rangle 
= \langle T \wh{f}, T \La^\dagger y \rangle
= \langle T \La^\dagger y - T [T_\cN^* T_\cN]^{-1} T_\cN^* T \La^\dagger y  , T \La^\dagger y \rangle \\
& = \|T \La^\dagger y \|^2 - \langle  [T_\cN^* T_\cN]^{-1} T_\cN^* T \La^\dagger y  , T_\cN^* T \La^\dagger y \rangle
= \|T \La^\dagger y \|^2 - ( T_\cN^* T \La^\dagger y )^*  [T_\cN^* T_\cN]^{-1} ( T_\cN^* T \La^\dagger y ).
\end{align*}
It follows that the squared Chebyshev center can be written as
\be
\label{Put1}
\rho^2 = \inf_{a \ge 0} \big\{  a (1-\|T \La^\dagger y \|^2)
+ ( a T_\cN^* T \La^\dagger y )^*  [a T_\cN^* T_\cN]^{-1} ( a T_\cN^* T \La^\dagger y )
: \;  a T_\cN^* T_\cN \succeq  \wt{Q}_\cN^* \wt{Q}_\cN  \big\}.
\ee
Relying on a classical property of Schur complements, we now notice that
\begin{align}
\label{Put2}
( a T_\cN^* T \La^\dagger y )^*  [a T_\cN^* T_\cN]^{-1} ( a T_\cN^* T \La^\dagger y )
& = \inf_{b \ge 0} \big\{ b : \; 
b \ge ( a T_\cN^* T \La^\dagger y )^*  [a T_\cN^* T_\cN]^{-1} ( a T_\cN^* T \La^\dagger y ) \big\}\\
& = 
\nonumber
\inf_{b \ge 0} \left\{ b : \;
\bbmx
a T_\cN^* T_\cN & \vline & a T_\cN^* T \La^\dagger y\\
\hline 
( a T_\cN^* T \La^\dagger y )^* & \vline & b
\ebmx \succeq 0
\right\}.
\end{align}
Putting \eqref{Put1} and \eqref{Put2} together,
we now see that $\rho^2$ has the semidefinite expression
$$
\rho^2
= \inf_{a,b \ge 0} \left\{  
a (1-\|T \La^\dagger y \|^2) + b : \;
a T_\cN^* T_\cN \succeq  \wt{Q}_\cN^* \wt{Q}_\cN,
\;
\bbmx
a T_\cN^* T_\cN & \vline & a T_\cN^* T \La^\dagger y\\
\hline 
( a T_\cN^* T \La^\dagger y )^* & \vline & b
\ebmx \succeq 0
\right\}.
$$
Next,
for the two-hyperellipsoid-intersection,
our upper bound for the Chebyshev radius is the minimal value of $ \la(\tau) \times  \big( 1 - ( (1-\tau) \|R f^\tau\|^2 + \tau \| S f^\tau \|^2 ) \big)$ over $\tau \in [0,1]$.
For a fixed $\tau$,
this reduces to the previous situation by considering $T = [(1-\tau) R^* R + \tau S^* S]^{1/2}$,
yielding a semidefinite expression with variables $a,b \ge 0$ that ought to be minimized over $\tau \in [0,1]$ as well.
By making the change of variables $c = a(1-\tau)$
and $d = a \tau$,
we arrive, after some work,
at the semidefinite program announced in \eqref{SDP2Hyper}.

\item 
As noted in a),
the optimal parameter $\tau^y$ can be found by solving the semidefinite program~\eqref{SDP2Hyper} before setting $\tau^y = d^y / (c^y + d^y)$.
This parameter can alternatively be obtained by solving the equation $d\|\wt{Q} h^\tau\|^2 / d\tau = 0$---recall that $\|\wt{Q} h^\tau\|^2 = (1-G(\tau))/H(\tau)$ and (see \eqref{Der=0}) that $d\|\wt{Q} h^\tau\|^2 / d\tau  = F(\tau)/H(\tau)$,
where $F(\tau) := \big( \|R f^\tau\|^2 + \|R h^\tau\|^2 \big) - \big( \|S f^\tau\|^2 + \|S h^\tau\|^2 \big)$.
Thus, we need to solve the equation $F(\tau ) = 0$,
which can be done via the Newton iteration
$$
\tau \leftarrow \tau - \f{F(\tau)}{F'(\tau)}
$$
provided that $F(\tau)$ and $F'(\tau)$ are accessible throughout the iterations.
Note that performing an eigendecomposition at each iteration gives $H(\tau)= 1/\la(\tau)$ and $h^\tau$
and that $f^\tau$ can also be computed,
thus making $F(\tau)$ accessible.
As for $F'(\tau)$,
it can also become accessible in some cases,
e.g. when $\wt{Q} = \Id_F$.
To see this,
notice that
$$
F'(\tau) = 2 (\alpha - \beta),
\qquad
\alpha := \left\langle (R^* R - S^* S)f^\tau, \f{d f^\tau}{d \tau} \right\rangle,
\quad 
\beta := \left\langle (R^* R - S^* S)h^\tau, \f{d h^\tau}{d \tau} \right\rangle.
$$
We can access $\alpha$ via the formula $f^\tau = \La^\dagger y - \big[ (1-\tau) R_\cN^* R_\cN + \tau S_\cN^* S_\cN \big]^{-1} \big( (1-\tau) R_\cN^* R + \tau S_\cN^* S \big) \La^\dagger y$ written as 
$\big[ (1-\tau) R_\cN^* R_\cN + \tau S_\cN^* S_\cN \big] (f^\tau - \La^\dagger y) = - \big( (1-\tau) R_\cN^* R + \tau S_\cN^* S \big) \La^\dagger y$,
which differentiates to give,
after some work,
$$
\f{d f^\tau}{d \tau} = \big[ (1-\tau) R_\cN^* R_\cN + \tau S_\cN^* S_\cN \big]^{-1} (R_\cN^* R - S_\cN^* S) f^\tau.
$$
Turning to $\beta$, we notice that the eigenequation 
$\big[ (1-\tau) R_\cN^* R_\cN + \tau S_\cN^* S_\cN \big] h^\tau = H(\tau) \wt{Q}_\cN^* \wt{Q}_\cN h^\tau$
yields
$$
\left\langle h^\tau, \f{d h^\tau}{d \tau} \right\rangle - \tau \left\langle (R_\cN^* R_\cN - S_\cN^* S_\cN ) h^\tau, \f{d h^\tau}{d \tau} \right\rangle
= H(\tau) \left\langle \wt{Q}_\cN^* \wt{Q}_\cN h^\tau, \f{d h^\tau}{d\tau} \right\rangle,
$$
while 
$$
\left\langle \wt{Q}_\cN^* \wt{Q}_\cN h^\tau, \f{d h^\tau}{d\tau} \right\rangle = \f{1}{2} \f{d \|\wt{Q}h^\tau \|^2}{d \tau} = \f{F(\tau)}{2 H(\tau)},
$$
so that $\beta$ is accessible as soon as $\left\langle h^\tau, \df{d h^\tau}{d \tau} \right\rangle$ is.
But when $\wt{Q} = \Id_F$, the latter is indeed accessible since
$$
\left\langle h^\tau, \f{d h^\tau}{d \tau} \right\rangle
= \f{1}{2} \f{d \| h^\tau\|^2}{d\tau}
= \f{1}{2} \f{d \|\wt{Q} h^\tau\|^2}{d\tau}
= \f{F(\tau)}{2 H(\tau)}.
$$

\item
It took considerable effort to determine the candidate locally optimal recovery map
and, even then,
we need an orthogonality condition to guarantee that we have found the genuine locally optimal recovery map.
In contrast,
it would take much less effort determine a locally near-optimal recovery map.
Indeed, any map assigning to $y \in \bR^m$ a model- and data-consistent $f^y \in F$---i.e., such that $\|Rf^y\| \le 1$, $\|S f^y \| \le 1$, and $\La f^y = y$---has a local worst-case error that is away from optimality by a factor~two. 
To see this, with $\wh{f}^y$ denoting the genuine Chebyshev center, we notice that 
\begin{align*}
{\rm lwce}_y(f^y) 
& = \sup_{\substack{\|Rf\| \le 1, \|S f\| \le 1 \\ \La f = y}} \|Q(f) - Q(f^y) \|
\le \sup_{\substack{\|Rf\| \le 1, \|S f\| \le 1 \\ \La f = y}} \|Q(f) - Q(\wh{f}^y)\|  \; + \;  \|Q(\wh{f}^y) - Q(f^y) \| \\
& \le 2 \sup_{\substack{\|Rf\| \le 1, \|S f\| \le 1 \\ \La f = y}} \|Q(f) - Q(\wh{f}^y)\|
 = 2 \, {\rm lwce}_y(\wh{f}^y).
\end{align*}
Such a model- and data-consistent $f^y$ can be obtained by solving the quadratic program
$$
\minimize{f \in F} \max\big\{ \|Rf\|, \|Sf\| \big\}
\qquad \mbox{s.to} \quad \La f = y.
$$
Using arguments similar to the ones from \cite[Theorem 4]{FouLiaVel},
we can show that the solution $f^y$ to the above program actually agrees with a constrained regularizer $f^\tau$ for some $\tau \in [0,1]$.

\end{enumerate}

\erk

\section{Locally Optimal Recovery in Hilbert Spaces}
\label{SecLoc}

This section finally deals with the third multi-fidelity  scenario mentioned in Section \ref{SecForma}.
Still working in Hilbert spaces
with $f_0$ and $f_0 - f_1$ once more modeled as belonging to the hyperellipsoids $\cK_0$ and $\cK_1$ from \eqref{2HyperEll},
we now target locally optimal recovery maps,
which are more challenging to produce than the globally optimal recovery maps treated in Section \ref{SecGlo}.
More precisely, given a linear quantity of interest $Q: F \to Z$,
we would like to construct,
for every $y=[y_0;y_1] \in \bR^{m_0+m_1}$,
an element $\Delta^{\rm opt}(y) = z^{\rm opt}$ that minimizes over all $z \in Z$ the local worst-case error 
$$
{\rm lwce}_y(z)
= \sup_{\substack{\|P_0 f_0\| \le \eps_0, \|P_1(f_0 - f_1)\| \le \eps_1\\ \La_0 f_0 = y_0 , \La_1 f_1 = y_1}} \|Q(f_0) - z\|.
$$
We capitalize on the result from the previous section by reformulating the latter as
\be
\label{ReformLWCE}
{\rm lwce}_y(z)
= \sup_{\substack{\|Rf\| \le 1, \|Sf\| \le 1\\ \La f = y}} \|\wt{Q}(f) - z\|,
\ee
where the linear maps $R$, $S$, $\La$, and $\wt{Q}$ are the ones introduced in \eqref{DefRS}.
Our main result shows that, under some specific assumptions, constrained regularization yet again provides an optimal recovery map,
but this time the optimal parameter $\tau^y \in [0,1]$ depends on $y$,
making the locally optimal recovery map nonlinear.

\bthm
\label{ThmLocOpt}
Suppose that $F$ is a Hilbert space and that the model sets $\cK_0,\cK_1 \inc F$
are the hyperellipsoidal sets defined in \eqref{2HyperEll}.
Suppose also that $\ker(\La_0) \inc \ker(P_1)$ or $\ker(\La_0) \inc \ker(\La_1)$.
If the linear quantity of interest $Q: F \to Z$ maps into a Hilbert space, then a locally optimal recovery map is given by
$\Delta^{\rm opt}: y \in \bR^{m_0+m_1}
\mapsto Q(f^{\tau^y}_0) \in Z$,
where
$$
[f^{\tau}_0;f^{\tau}_1]
:= \underset{f = [f_0; f_1] \in F \times F}{\argmin \;}
\left[ (1-\tau) \|P_0 f_0\|^2 + \tau \|P_1(f_0 - f_1) \|^2 \right]
\qquad \mbox{s.to} \quad \La_0 f_0 = y_0, \; \La_1 f_1 = y_1.
$$
Here, the regularization parameter $\tau^y$ is $\tau^y = c_1^y/(c_0^y+c_1^y)$,
with $b^y,c_0^y,c_1^y \ge 0$ being solutions to the following semidefinite program featuring $\cN_0 := \ker(\La_0)$, $\cN_1 := \ker(\La_1)$,
$u_0 := P_{0} \La_0^\dagger y_0$, 
and $u_1 := P_{1} (\La_0^\dagger y_0 - \La_1^\dagger y_1)$:
\begin{align}
\label{SDPforLocMult}
& \minimize{b,c_0,c_1 \ge 0} \; 
c_0 \big( \eps_0^2 - \|u_0\|^2 \big) + c_1 \big( \eps_1^2 - \|u_1\|^2 \big) + b\\
\nonumber
\mbox{s.to } &
\bbmx
c_0 P_{0,\cN_0}^* P_{0,\cN_0} + c_1 P_{1,\cN_1}^* P_{1,\cN_1} & \vline & - c_1 P_{1,\cN_1}^* P_{1,\cN_1}\\
\hline
- c_1 P_{1,\cN_1}^* P_{1,\cN_1} & \vline & c_1 P_{1,\cN_1}^* P_{1,\cN_1}
\ebmx \succeq 
\bbmx
Q_{\cN_0}^* Q_{\cN_0} & \vline & 0 \; \;\\
\hline 0 & \vline & 0 \; \;
\ebmx
\\
\nonumber
\mbox{and } &
\bbmx
c_0 P_{0,\cN_0}^* P_{0,\cN_0} + c_1 P_{1,\cN_1}^* P_{1,\cN_1} & \vline & - c_1 P_{1,\cN_1}^* P_{1,\cN_1}
& \vline & 
c_0 P_{0,\cN_0}^* u_0 + c_1 P_{1,\cN_1}^* u_1
\\
\hline
- c_1 P_{1,\cN_1}^* P_{1,\cN_1} & \vline & c_1 P_{1,\cN_1}^* P_{1,\cN_1} & \vline & - c_1 P_{1,\cN_1}^* u_1\\
\hline 
(c_0 P_{0,\cN_0}^* u_0 + c_1 P_{1,\cN_1}^* u_1)^* & \vline & - (c_1 P_{1,\cN_1}^* u_1)^* & \vline & b
\ebmx \succeq 
0. 
\end{align}
\ethm

The key is to leverage the additional assumptions to certify that the necessary (and sufficient) orthogonality condition from Theorem \ref{ThmCCtwoHyp} is met,
as shown below.

\blem
\label{PropOrthCond}
If $\ker(\La_0) \inc \ker(P_1)$ or $\ker(\La_0) \inc \ker(\La_1)$, then,
 for any $\tau \in (0,1)$,
$$
\langle P_0 f_0^\tau, P_0 h_0 \rangle  = 0
\quad \mbox{and} \quad
\langle P_1(f_0^\tau-f_1^\tau), P_1(h_0 - h_1) \rangle = 0
\quad \mbox{for all }  h_0 \in \ker(\La_0), \;  h_1 \in \ker(\La_1).
$$
\elem

\bpf
With the notation of \eqref{DefRS}, 
we see that $f^\tau = [f_0^\tau;f_1^\tau]$ is the minimizer of $(1-\tau) \|Rf\|^2 + \tau \|Sf\|^2$ subject to $\La f = y$.
As such, it must satisfy 
$\big( (1-\tau) R^* R + \tau S^* S \big) f^\tau \perp \ker(\La)$,
in other words,  $(1-\tau) \langle R f^\tau, R h \rangle + \tau \langle S f^\tau, S h \rangle = 0 $ for all $h = [h_0; h_1] \in \ker \La$.
Reinserting the specific expressions of $R$ and $S$,
this reads
$$
\f{1-\tau}{\eps_0^2} \langle P_0(f_0^\tau), P_0(h_0) \rangle 
+ \f{\tau}{\eps_1^2} \langle P_1(f_0^\tau-f_1^\tau), P_1(h_0 - h_1) \rangle =0 
\qquad \mbox{for all } h_0 \in \ker ( \La_0), \, h_1 \in \ker( \La_1 ),
$$ 
which itself decouples as the two identities
\begin{align}
\label{Iden1}
\f{1-\tau}{\eps_0^2}  \langle P_0 f_0^\tau, P_0 h_0 \rangle + \f{\tau}{\eps_1^2} \langle P_1 (f_0^\tau - f_1^\tau), P_1 h_0 \rangle  & = 0 
\qquad \mbox{for all } h_0 \in \ker \La_0,\\
\label{Iden2}
\f{\tau}{\eps_1^2} \langle  P_1 (f_0^\tau - f_1^\tau), P_1 h_1 \rangle & = 0
\qquad \mbox{for all } h_1 \in \ker \La_1.
\end{align}
First, let us assume that $\ker(\La_0) \inc \ker(P_1)$.
Then \eqref{Iden1} directly implies that $\langle P_0 f_0^\tau, P_0 h_0 \rangle = 0$ for all $h_0 \in \ker(\La_0)$,
while $\langle P_1(f_0^\tau-f_1^\tau), P_1(h_0 - h_1) \rangle = - \langle P_1(f_0^\tau-f_1^\tau), P_1  h_1 \rangle = 0$ for all $h_0 \in \ker(\La_0)$ and $h_1 \in \ker(\La_1)$ according to \eqref{Iden2}.
Second, let us now assume that $\ker(\La_0) \inc \ker(\La_1)$.
Then \eqref{Iden2} implies that $\langle P_1 (f_0^\tau - f_1^\tau), P_1 h_0 \rangle = 0$, which,
when substituted in \eqref{Iden1},
yields $\langle P_0 f_0^\tau, P_0 h_0 \rangle = 0$ for all  $h_0 \in \ker(\La_0)$
and, when combined with \eqref{Iden2},
yields $\langle P_1 (f_0^\tau - f_1^\tau), P_1(h_0 - h_1) \rangle = 0$ for all $h_0 \in \ker(\La_0)$ and $h_1 \in \ker(\La_1)$.
\epf

We are now ready to justify the main result of this section.

\bpf[Proof of Theorem \ref{ThmLocOpt}]
For $y \in \bR^{m_0+m_1}$,
in view of the reformulation \eqref{ReformLWCE} of the local worst-case error,
the desired $\Delta^{\rm opt}(y)$ is the Chebyshev center of the set $\wt{Q}(\{ f \in F: \|R f \| \le 1, \|S f\| \le 1, \La f = y \})$.
By Theorem \ref{ThmCCtwoHyp} and using the notation introduced there,
the latter is equal to $\wt{Q}({f^{\tau^y}})$,
i.e., $Q(f_0^{\tau^y})$, provided that $\langle R f^{\tau^y}, R h^{\tau^y} \rangle = 0$,
i.e., $\langle P_0 f_0^{\tau^y}, P_0 h_0^{\tau^y} \rangle = 0$.
But Lemma \ref{PropOrthCond} guarantees that this orthogonality condition holds for any $h_0 \in \ker(\La_0)$ and in particular for $h_0^{\tau^y}$.
Now, according to a) in the remark closing Section \ref{SecCCtwoHyp},
the parameter $\tau^y$ is obtained as $\tau^y = c_1^y / (c_0^y + c_1^y)$ where $b^y,c_0^y,c_1^y \ge 0$ solve \eqref{SDP2Hyper}.
By substituting the expressions for $R$, $S$, $\La$, and $\wt{Q}$ from \eqref{DefRS},
we transform this semidefinite program,
after some routine work,
into the announced semidefinite program \eqref{SDPforLocMult}.
\epf

\end{document}